\documentclass{article}

\newtheorem{theorem}{Theorem}%[section]
\newtheorem{lemma}[theorem]{Lemma}
\newtheorem{proposition}[theorem]{Proposition}
\newtheorem{corollary}[theorem]{Corollary}

\newtheorem{remark}[theorem]{Remark}
\newtheorem{example}[theorem]{Example}
\newtheorem{notation}[theorem]{Notation}

\newenvironment{proof}[1][Proof:]{\begin{trivlist}
\item[\hskip \labelsep {\bfseries #1}]}{\end{trivlist}}

\newcommand{\qed}{\nobreak \ifvmode \relax \else
      \ifdim\lastskip<1.5em \hskip-\lastskip
      \hskip1.5em plus0em minus0.5em \fi \nobreak
      \vrule height0.75em width0.5em depth0.25em\fi}

 0 3 \font\ccc =msbm10

\begin{document}

\title{Semi-symmetric algebras: General Constructions. Part II}

\author{Valentin Vankov Iliev\\
Section of Algebra,\\ 
Institute of Mathematics and Informatics,\\
1113 Sofia, Bulgaria}

\maketitle

\begin{abstract}

In \cite{[3]} we present the construction of the semi-symmetric
algebra $[\chi](E)$ of a module $E$ over a commutative  ring  $K$
with  unit, which generalizes the tensor algebra $T(E)$, the
symmetric algebra $S(E)$, and the exterior algebra $\wedge(E)$,
deduce some of its functorial properties, and prove a
classification theorem. In the present paper we continue the study
of the semi-symmetric algebra and discuss its graded dual, the
corresponding canonical bilinear form, its coalgebra structure, as
well as left and right inner products. Here we present a unified
treatment of these topics whose exposition in \cite[A.III]{[2]} is
made simultaneously for the above three particular (and, without a
shadow of doubt --- most important) cases.

\end{abstract}

\section{Introducton}

\label{I}

In order to make  the  exposition  self-contained, in this
introduction we remind the main definitions and
results from \cite{[3]}.

Let $K$ be a commutative ring with unit $1$. Denote by $U(K)$ the
group of units of $K$. Given a positive integer $d$, let $W\leq
S_d$ be a permutation group, and let $\chi$ be a linear $K$-valued
character of the group $W$, that is, a group homomorphism
$\chi\colon W\to U(K)$. We call a $W$-module any $K$-linear
representation of $W$  and view it also as a left unitary module
over the group ring $KW$. Let $M$ be a $W$-module. We denote by
${}_\chi M$ the $W$-submodule of $M$, generated by all differences
$\chi(\sigma)z-\sigma z$, where $\sigma\in W$, $z\in M$, and by
$M_\chi$ the $W$-submodule of $M$, consisting of all $z\in M$ such
that $\sigma z=\chi(\sigma)z$ for all $\sigma\in W$. Given
$K$-modules $E$, $F$, we denote  by  $Mult_K(E^d,F)$ the
$K$-module consisting of all $K$-multilinear maps $E\to F$, and by
$T^d(E)$ --- the $d$-th tensor power of $E$. The $K$-modules
$T^d(E)$, $Hom_K(T^d(E),F)$, and  $Mult_K(E^d,F)$ have the  usual
structure  of $W$-modules, see \cite[Ch. III, Sec. 5, $n^o$
1]{[1]}. We denote the factor- module
$T^d(E)/{}_{\chi^{-1}}T^d(E)$ by $[\chi]^d(E)$, and call it
\emph{$d$-th semi-symmetric power of weight $\chi$ of the
$K$-module $E$}. By definition, $[\chi]^0(E)=K$. The image of the
tensor $x_1\hskip-2pt\otimes\ldots \otimes\hskip-1pt x_d\in
T^d(E)$ by the canonical  homomorphism $\varphi_d\colon T^d(E)\to
[\chi]^d(E)$ is denoted by $x_1\chi\ldots \chi x_d$, and is called
\emph{decomposable $d-\chi$-vector}. Thus,
$x_{\sigma\left(1\right)}\chi\ldots \chi x_{\sigma\left(d\right)}=
\chi(\sigma)x_1\chi\ldots \chi x_d$ for any permutation $\sigma\in
W$.

In \cite[(1.1.1)]{[3]} we show that $d$-th semi-symmetric power
$[\chi]^d(E)$ is a representing  object for the functor
$Mult_K(E^d,-)_\chi$. As usual, we denote by $S_\infty$ the group
of all permutations of the set of all positive integers, which fix
all but finitely many elements. We identify the symmetric group
$S_d$ with  the  subgroup  of $S_\infty$, consisting  of  all
permutations fixing any $n>d$.  Let $(W_d)_{d\geq 1}$ be a
sequence of subgroups of $S_\infty$. This sequence is said to be
\emph{admissible} if  $W_d\leq S_d$ for all $d\geq 1$. A sequence
of $K$-valued characters $(\chi_d\colon W_d\to U(K))_{d\geq 1}$ is
said to be \emph{admissible} if its sequence of domains
$(W_d)_{d\geq 1}$ is admissible. We define an injective
endomorphism $\omega$ of the symmetric group $S_\infty$ by the
formula $(\omega(\sigma))(d) =\sigma(d-1)+1$, $(\omega(\sigma))(1)=1$. A
sequence $(W_d)_{d\geq 1}$ is called \emph{$\omega$-stable} if it
is admissible, and $W_d\leq W_{d+1}$, $\omega(W_d)\leq W_{d+1}$,
for all $d\geq 1$. A  sequence  of linear $K$-valued characters
$(\chi_d\colon W_d\to U(K))_{d\geq 1}$ is said to be
\emph{$\omega$-invariant} if its sequence of domains $(W_d)_{d\geq
1}$ is $\omega$-stable, and
\[
{\chi_{d+1}}_{\mid W_d}=\chi_d={\chi_{d+1}}\circ\omega_{\mid W_d}
\]
for all $d\geq 1$. Given  a  $K$-module $E$, any admissible
sequence of characters  $\chi=(\chi_d)_{d\geq 1}$ produces a
graded $K$-module $[\chi](E)=\coprod_{d\geq 0}[\chi]^d(E)$, where
$[\chi]^d(E)=[\chi_d]^d(E)$, and $[\chi]^0(E)=K$. Denote by
$\varphi(E)$ the  canonical  $K$-linear  homomorphism
$\coprod_{d\geq 0}\varphi_d\colon T(E)\to [\chi](E)$, where
$\varphi_0=id_K$. We denote by $K^{\left(\infty\right)}$ a free
$K$-module with countable basis. The following two theorems are
proved in \cite{[3]} (see \cite[1.3.1]{[3]} and
\cite[1.3.3]{[3]}):

\begin{theorem}\label{I.1}  Let $\chi$ be an admissible  
sequence  of  characters.  The
following statements are then equivalent.

(i) The sequence $\chi$ is $\omega$-invariant;

(ii) for any $K$-module $E$ the $K$-module $[\chi](E)$  has  a  structure  of
associative graded $K$-algebra, such that $\varphi(E)$ is a homomorphism  of
graded $K$-algebras;

(iii) the $K$-module $[\chi](K^{\left(\infty\right)})$ has a
structure of associative graded $K$-algebra, such that
$\varphi(K^{\left(\infty\right)})$ is  a  homomorphism of graded
$K$-algebras.

\end{theorem}

The $K$-algebra $[\chi](E)$  is  called  the \emph{semi-symmetric
algebra of weight $\chi$  of the $K$-module $E$}, and its
elements --- \emph{$\chi$-vectors}.

\begin{theorem}\label{I.2} Let $W=(W_d)_{d\geq 1}$ be an $\omega$-stable
sequence of groups. Then the group of all $\omega$-invariant sequences of
characters on $W$ (with componentwise multiplication) is trivial  or
isomorphic to the multiplicative subgroup of $K$ consisting of all involutions.

\end{theorem}

We obtain immediately

\begin{corollary} \label{I.3} If $\chi=(\chi_d)_{d\geq 1}$
is an $\omega$-stable sequence of characters, then

(i) one has $\chi=\chi^{-1}$, where
$\chi^{-1}=(\chi_d^{-1})_{d\geq 1}$;

(ii) if the ring $K$ is an integral domain, then the possible
values of $\chi_d$ in $K$ are $\pm 1$ for any $d\geq 1$.

\end{corollary}

When $W_d=\{1\}$ for all $d\geq 1$, the graded algebra $[\chi](E)$
coincides with the tensor algebra $T(E)$. When $W_d=S_d$  and
$\chi_d$ is the unit character for all $d\geq 1$, the graded
algebra $[\chi](E)$ coincides with the symmetric algebra $S(E)$.
When $W_d=S_d$ and $\chi_d$ is the signature for all $d\geq 1$,
the graded algebra $[\chi](E)$ is the anti-symmetric algebra of
$E$; in particular, if $1/2\in K$, then $[\chi](E)$ is the
exterior algebra $\wedge(E)$ of the $K$-module $E$. If $E$ is a
$n$-generated $K$-module, $k\geq n$, and if $W_d=\{1\}$ for all
$d\leq k$, $W_d=S_d$ for all $d>k$, and $\chi_d$ is the signature
for all $d\geq 1$, then $[\chi](E)$ is \emph{the tensor algebra
truncated by its elements of degree $> k$}.

Let $W\leq S_d$ be a permutation group and let $\chi$ be a linear
$K$-valued character of the group $W$. In \cite[1]{[5]} we
construct a basis for the $d$-th semi-symmetric power
$[\chi]^d(E)$, $d\geq 1$,  starting  from  the  standard basis for
$T^d(E)$ in the case $K$ is a field of characteristics $0$, but
the results hold when $K$ is a commutative ring with unit, which
is an integral domain, the order of the group $W$ is invertible in
$K$, and the $K$-module $E$ is free, see \cite{[4]} where this
generalization was announced. The counterexamples from \cite{[4]}
show that these conditions are necessary for $[\chi]^d(E)$ to be a
free $K$-module for all permutation groups $W\leq S_d$ and for all
characters $\chi\colon W\to U(K)$. Here we prove these general
results, see Theorem~\ref{II.6}, its Corollary~\ref{II.7}, and
Example~\ref{II.11}.

In this paper we continue the study of semi-symmetric algebras
under the condition that the commutative ring $K$ is both a
$\hbox{\ccc Q}$-ring and an integral domain, and under the
assumption that the $K$-module $E$ is a free $K$-module with a
finite basis. We unite the bases for $[\chi]^d(E)$, $d\geq 0$, and
get a basis for the semi-symmetric algebra $[\chi](E)$ of weight
$\chi$, considered as a $K$-module. This is done in
Corollary~\ref{II.10}.

Further, we study some duality properties of the semi-symmetric
powers and algebras of weight $\chi$. In Theorem~\ref{III.1} we
define a non-singular bilinear form on the product
$[\chi]^d(E)\times [\chi^{-1}]^d(E^*)$, and use it to identify the
K-modules $([\chi]^d(E))^*$ and $[\chi^{-1}](E^*)$. Mimicking the
case of an exterior power, we make use of generalized Schur
function (see \cite{[6]}) instead of determinant. After this
identification, the above bilinear form coincides with the
canonical bilinear form of the $K$-module $[\chi]^d(E)$; here
$M^*$ denotes the dual of the $K$-module $M$. Thus, we get an
identification of the semi-symmetric algebra $[\chi](E^*)$ with
the dual graded algebra $([\chi](E))^{*gr}$ of  the semi-symmetric
algebra $[\chi](E)$, see Theorem~\ref{III.13}, (i). Moreover, we
extend the sequence of the above canonical bilinear forms to the
canonical bilinear form of the graded algebras $[\chi](E)$ and
$([\chi](E))^{\otimes k}$, by assuming that the homogeneous
components are orthogonal, see Theorem~\ref{III.13}, (ii), (iii).
Because of the above identification, the elements of the
semi-symmetric algebra $[\chi](E^*)$ are called
\emph{$\chi$-forms}. In Corollary~\ref{IV.2} we define a structure
of graded coassociative and counital $K$-coalgebra on $[\chi](E)$,
and show that the structure of graded associative algebra with
unit on its dual $([\chi](E))^{*gr}=[\chi](E^*)$, defined by
functoriality, coincide with the usual structure of graded
associative algebra with unit on the graded $K$-module
$[\chi](E^*)$. In particular, when $[\chi](E)$ is the graded
$K$-module underlaying the symmetric algebra (or  the exterior
algebra, or the tensor algebra) of the $K$-module $E$, we obtain
the usual structure of $K$-coalgebra on it (see~\cite[A III,
139-141]{[2]}). In Section~\ref{V}, following \cite[Ch. III, Sec.
8, $n^o$ 4]{[1]}, we find out the main properties of the left and
right inner products of a $\chi$-vector and a $\chi$-form.

\section{Basis of semi-symmetric algebra\\
 of a free module}

\label{II}

Let $W$ be a  a finite group, and let $\chi$ be a
linear $K$-valued character of the group $W$.
Let us assume that $|W|\in U(K)$ and set $a_\chi=|W|^{-1}\sum_{\sigma\in W}\chi^{-1}
(\sigma)\sigma$. The element $a_\chi$ of the group ring
$KW$ defines $K$-linear endomorphism $a_\chi\colon M\to
M$ by the rule $z\mapsto a_\chi z$. Then
the $W$-submodule ${}_\chi M$ of $M$ is the kernel of $a_\chi$, and
the $W$-submodule $M_\chi$ of $M$ is the image of $a_\chi$.

Let $M$ be a free $K$-module with basis $(e_i)_{i\in I}$.
Let us suppose that the finite group $W$ acts on the index set $I$.
Denote by $W_i$ the stabilizer of $i\in I$ and by $W^{\left(i\right)}$
a system of representatives of the left classes of $W$ modulo $W_i$.
Let $(\gamma_i)_{i\in I}$ be a family of maps
$W\to U(K)$ such that
$\gamma_i(\sigma\tau)=\gamma_{\tau i}(\sigma)\gamma_i(\tau)$
for all $i\in I$, and all $\sigma, \tau\in W$.
In particular, the restriction of $\gamma_i$ on $W_i$ is a
$K$-valued character of the group $W_i$ for any $i\in I$.
The $K$-module $M$ has a structure of monomial $W$-module, defined by the rule
\begin{equation}
\sigma e_i=\gamma_i(\sigma)e_{\sigma i},\mbox{\ }\sigma\in W,
\mbox{\ }i\in I. \label{II.1}
\end{equation}

We set $I(\chi,M)=\{i\in I\mid \gamma_i=\chi\mbox{\  on\ } W_i\}$,
$I_0(\chi,M)=I\backslash I(\chi,M)$.

\begin{lemma} \label{II.2} (i) The set $I(\chi,M)$ is a
$W$-stable subset of $I$;

(ii) one has $a_\chi (v_i)=0$ for $i\in I_0(\chi,M)$.

\end{lemma}

\begin{proof} (i) Given $i\in I$, suppose $\sigma\in W$ and
$\tau\in W_i.$ Then $W_{\sigma i}=\sigma W_i\sigma^{-1}$ and $\chi
(\sigma\tau\sigma^{-1})=\chi (\tau)$. Moreover,
\[
\gamma_{\sigma i}(\sigma\tau\sigma^{-1})= \gamma_{\sigma^{-1}\sigma
i}(\sigma\tau)\gamma_{\sigma i}(\sigma^{-1})= \gamma_{\sigma
i}(\sigma^{-1})\gamma_i(\sigma\tau)=
\]
\[
\gamma_{\sigma\tau i}(\sigma^{-1})\gamma_i(\sigma\tau)=
\gamma_i(\sigma^{-1}\sigma\tau)= \gamma_i(\tau).
\]

(ii) The complement of $I(\chi,M)$ in $I$ also is $W$-stable; let
$i\in I\backslash I(\chi,M)$. We have
\[
a_\chi (v_i)=|W|^{-1}\sum_{\sigma\in
W^{(i)}}\sum_{\tau\in W_i}
\chi^{-1}(\sigma\tau)\gamma_i(\sigma\tau)v_{\sigma\tau i}=
\]
\[
|W|^{-1}\sum_{\sigma\in W^{(i)}}\chi^{-1}
(\sigma)\gamma_i(\sigma)(\sum_{\tau\in
W_i}\chi^{-1}(\tau)\gamma_i(\tau))v_{\sigma i},
\]
and the equality $a_\chi (v_i)=0$ holds because the product
$\chi^{-1}\gamma_i$ is not the unit character of the group $W_i$.

\end{proof}

We choose an element $i$ from any $W$-orbit in $I$ and denote the
set of these $i$'s by $I^*$. Finally,  we set $J(\chi,M)=I^*\cap
I(\chi,M)$, and $J_0(\chi,M)=I^*\cap I_0(\chi,M)$. Following
\cite[Ch. III, Sec. 5, $n^o$ 4]{[2]}, we get a basis of the
$K$-module $M$ consisting of
\begin{eqnarray}
\mbox{ } & \mbox{ } & e_j,\quad j\in J(\chi,M),\label{II.3} \\
\mbox{ } & \mbox{ } & e_i-\chi(\sigma)\gamma_i(\sigma)e_{\sigma i},\quad i\in
I^*,\ \sigma\in W^{(i)}, \ \sigma\not\in W_i, \label{II.4} \\
\mbox{ } & \mbox{ } & e_i,\quad i\in J_0(\chi,M).\label{II.5}
\end{eqnarray}

\begin{theorem}\label{II.6} Let the ring $K$ be an integral domain and let $|W|\in
U(K)$. Then

(i) the union of the families~(\ref{II.4})
and~(\ref{II.5}) is a basis for ${}_\chi M$;

(ii) the family $a_\chi (e_j)$, \ $j\in J(\chi,M)$,
is a basis for $M_\chi $;

(iii) the family $e_j\ mod({}_\chi M)$, \ $j\in J(\chi,M)$, is a
a basis for the factor-module $M/{}_\chi M$.

\end{theorem}

\begin{proof} (i) The family~(\ref{II.4}) is in ${}_\chi M$ by definition. Lemma~\ref{II.2},
(ii), implies that the family~(\ref{II.5}) is contained in
${}_\chi M.$ Now, set $J=J(\chi,M)$ and suppose that $\sum_{j\in
J}k_ja_\chi (v_j)=0$ for some $k_j\in K$ such that $k_j=0$ for all
but a finite number of indices $j\in J$. We have
$$
\sum_{j\in J}k_ja_\chi (v_j)=|W|^{-1}\sum_{j\in J}\sum_{\sigma\in
W^{(j)}}k_j|W_j|\chi^{-1}(\sigma)\gamma_j(\sigma)v_{\sigma j},
$$
hence $k_j=0$ for all $j\in J$, which proves part (i). In
addition, we have proved that the elements $a_\chi (v_j)$, $j\in
J(\chi,M)$, are linearly independent.

(ii) The elements $a_\chi (v_j)$, $j\in J(\chi,M)$, are in
$M_\chi$ and, moreover, each element of $M_\chi$ has the form
$a_\chi (z)$ for some $z\in M$. Since the union of
families~(\ref{II.3}) --~(\ref{II.5}) is a basis for $M$ and since
the endomorphism $a_\chi$ annihilates~(\ref{II.4})
and~(\ref{II.5}), part (ii) holds.

(iii) Part (ii) implies part (iii).

\end{proof} 

Now, let us suppose that the $K$-module $E$ has basis
$(e_\ell)_{\ell\in L}$. Then the tensor power $M=T^d(E)$ has basis
$(e_i)_{i\in L^d}$, and if $W\leq S_d$ is a permutation group, the
rule $\sigma e_i=e_{\sigma i}$, $\sigma\in W$, defines on $M$ a
structure of monomial $W$-module.

\begin{corollary}\label{II.7} Let $W\leq S_d$  be a permutation group and let $\chi$
be  a linear $K$-valued character of $W$. If $K$ is an integral domain and $|W|\in U(K)$,  then
the $d$-th semi-symmetric power $[\chi](E)$ of weight $\chi$
of a free $K$-module $E$ with basis $(e_\ell)_{\ell\in L}$
is a free $K$-module with basis
\[
(e_{j_1}\chi\ldots\chi e_{j_d})_{\left(j_1,\ldots,j_d\right)\in J\left(\chi,T^d\left(E\right)\right)}.
\]

\end{corollary}

\begin{proof} Substitute $M=T^d\left(E\right)$, $I=L^d$, $\gamma_i(\sigma)=1$ for all $\sigma\in
W$, $i\in L^d$, in Theorem~\ref{II.6}.

\end{proof}

\begin{corollary}\label{II.8} Let $W\leq S_d$  be a permutation group and let $\chi$
be  a linear $K$-valued character of $W$. If $K$ is an integral
domain, $|W|\in U(K)$, and if $E$ is a projective $K$-module (a
projective $K$-module of finite type), then the $d$-th
semi-symmetric power $[\chi](E)$ of weight $\chi$ is a projective
$K$-module (a projective $K$-module of finite type).

\end{corollary}

\begin{proof} Let $L$ be a set (a finite set), and let
$K^{\left(L\right)}$ be the free $K$-module with the canonical
basis indexed by $L$. Let
\[
0\to E\to K^{\left(L\right)}
\]
be a splitting monomorphism of $K$-modules. Since the functor
$[\chi]^d(-)$ transforms epimorphisms into epimorphisms, the
sequence
\[
0\to [\chi]^d(E)\to [\chi]^d(K^{\left(L\right)})
\]
also is a splitting monomorphism of $K$-modules, and, moreover,
according to Corollary~\ref{II.7}, $[\chi]^d(K^{\left(L\right)})$
is a free $K$-module (free module with finite basis). Therefore
$[\chi]^d(E)$ is a projective $K$-module (a projective $K$-module
of finite type).

\end{proof}

\begin{remark}\label{II.9} Let us set
$J\left(\chi,T^0\left(E\right)\right)=\{\emptyset\}$,
$e_\emptyset=1$. We unite the bases of all semi-symmetric powers
$[\chi]^d(E)$ (see Corollary~\ref{II.7}), thus getting
$J\left(\chi,T\left(E\right)\right)=\cup_{d\geq
0}J\left(\chi_d,T^d\left(E\right)\right)$. In particular, when
$L=[1,n]$, the  elements of the set
$J\left(\chi,T^d\left(E\right)\right)$ can be chosen to be
lexicographically minimal in their $W$-orbits, and we can
introduce following notation:
\[
I\left(\chi,T^d\left(E\right)\right)=I(\chi,n,d),\mbox{\ }
I_0\left(\chi,T^d\left(E\right)\right)=I_0(\chi,n,d),
\]
\[
J\left(T^d\left(E\right),\chi\right)=J(\chi,n,d), \mbox{\ }
J\left(T\left(E\right),\chi\right)=J(\chi,n).
\]
For any $i\in I(\chi,n,d)$ we define $\ell m(i)$ to be the
lexicographically minimal element in the $W$-orbit of $i$, and set
$\zeta(i)=\chi_d(\sigma)$, where $\sigma\in W_d$ is such that
$\sigma i=\ell m(i)$. Since the restriction of the character
$\chi_d$ is identically $1$ on the stabilizer $(W_d)_i$, the
element $\zeta(i)\in U(K)$ does not depend on the choice of
$\sigma$.

\end{remark}

Let $\chi=(\chi_d)_{d\geq 1}$ be an $\omega$-invariant sequence of
characters and  let $W=(W_d)_{d\geq 1}$ be the sequence of their
domains.

\begin{corollary}\label{II.10} Let $K$ be both a $\hbox{\ccc Q}$-ring
and an integral domain.

(i) If $E$ is a $K$-module with basis $(e_\ell)_{\ell\in L}$, then
the family $(e_j)_{j\in J\left(T\left(E\right),\chi\right)}$ is a
basis for the semi-symmetric algebra $[\chi](E)$ of weight $\chi$,
considered as a $K$-module;

(ii) If $E$ is a $K$-module with finite basis
$(e_\ell)_{\ell=1}^n$, then the family $(e_j)_{j\in
J\left(\chi,n\right)}$ is a basis for the semi-symmetric algebra
$[\chi](E)$ of weight $\chi$, considered as a $K$-module. If $j\in
J(\chi,n,d)$ and $k\in J(\chi,n,e)$, then the multiplication table
of the $K$-algebra $[\chi](E)$ is given by the formulae
\[
e_j\chi e_k=
\left\{
\begin{array}{ll}
0 & \mbox{if $(j,k)\in I_0(\chi,n,d+e)$}\\
\zeta(j,k)e_{\ell m\left(j,k\right)} & \mbox{if $(j,k)\in
I(\chi,n,d+e)$}.
 \end{array}
\right.
\]

\end{corollary}

\begin{proof} (i) Straightforward use of Corollary~\ref{II.7}.

(ii) The first part is a particular case of (i). We have $e_j\chi
e_k=e_{\left(j,k\right)}$, and in case $(j,k)\in I_0(\chi,n,d+e)$
Lemma~\ref{II.2}, (ii), implies $e_{\left(j,k\right)}=0$.
Otherwise, $e_{\ell m\left(j,k\right)}\in J(\chi,n,d+e)$, and we
make use of Remark~\ref{II.9}.

\end{proof}

\begin{example}\label{II.11} We will show that if some of the conditions
of Corollary~\ref{II.7} fail, then the $K$-module $[\chi]^d(E)$ is
not necessarily free.

(i) The ring $K$ is not an integral domain.

We set $K=\hbox{\ccc Z}_{15}$,
$W=\{(1),(12)(34),(13)(24),(14)(23)\}\leq S_4$ is the Klein four
group, $\chi((12)(34))=4$, $\chi((13)(24))=4$, $\chi((14)(23))=1$,
$E=Ke_1\coprod Ke_2$, $I=[1,2]^4$,
$e_i=e_{i_1}\hskip-2pt\otimes\ldots\otimes\hskip-2pt e_{i_4}$ for
$i=(i_1,\ldots, i_4)\in I$.
 We have $\chi=\chi^{-1}$. The $K$-module $T^4(E)_\chi$ is spanned by the elements
\[
a_\chi(e_{\left(1,1,1,1\right)})=10e_{\left(1,1,1,1\right)},
\]
\[
a_\chi(e_{\left(2,2,2,2\right)})=10e_{\left(2,2,2,2\right)},
\]
\[
a_\chi(e_{\left(1,1,2,2\right)})=5e_{\left(1,1,2,2\right)}+5e_{\left(2,2,1,1\right)},
\]
\[
a_\chi(e_{\left(1,2,1,2\right)})=5e_{\left(1,2,1,2\right)}+5e_{\left(2,1,2,1\right)},
\]
\[
a_\chi(e_{\left(1,2,2,1\right)})=8e_{\left(1,2,2,1\right)}+8e_{\left(2,1,1,2\right)},
\]
\[
a_\chi(e_{\left(1,1,1,2\right)})=e_{\left(1,1,1,2\right)}+e_{\left(1,1,2,1\right)}+
e_{\left(1,2,1,1\right)}+4e_{\left(2,1,1,1\right)},
\]
\[
a_\chi(e_{\left(1,2,2,2\right)})=e_{\left(1,2,2,2\right)}+e_{\left(2,1,2,2\right)}+
e_{\left(2,2,1,2\right)}+4e_{\left(2,2,2,1\right)}.
\]
Thus, the $\hbox{\ccc Z}_{15}$-module $M_\chi$ is isomorphic to
the submodule
\[
\hbox{\ccc Z}_{15}e^{\left(1\right)}\coprod\hbox{\ccc
Z}_{15}e^{\left(2\right)}\coprod\hbox{\ccc
Z}_{15}e^{\left(3\right)}\coprod \hbox{\ccc
Z}_{15}10e^{\left(4\right)}\coprod\hbox{\ccc
Z}_{15}10e^{\left(5\right)}\coprod\hbox{\ccc
Z}_{15}5e^{\left(6\right)}\coprod\hbox{\ccc
Z}_{15}5e^{\left(7\right)}
\]
of a free $\hbox{\ccc Z}_{15}$-module with $7$ generators
$e^{\left(1\right)}$, $\ldots,$ $e^{\left(7\right)}$. This
submodule has $15^33^4$ elements, and this number is not a power
of $15$, hence $[\chi]^4(E)=M/{}_\chi M$ is not a free $\hbox{\ccc
Z}_{15}$-module.

(ii) The order $|W|$ of the group $W$ is not invertible in the
ring $K$.

We denote by $\varepsilon$ a primitive $3$-th root of unity and
set $K=\hbox{\ccc Z}[\varepsilon]$, $W=\{(1), (123), (132)\}\leq
S_3$, $\chi(123)=\varepsilon$, $E=Ke_1\coprod Ke_2$, $M=T^3(E)$,
$I=[1,2]^3$, $e_i=e_{i_1}\hskip-2pt\otimes\hskip-2pt
e_{i_2}\hskip-2pt\otimes\hskip-2pt e_{i_3}$ for $i=(i_1,i_2,
i_3)\in I$. The $K$-module $[\chi^2]^3(E)=M/{}_\chi M$ is spanned
by the elements
\[
e_{\left(1,1,1\right)}, e_{\left(2,2,2\right)},
e_{\left(1,1,2\right)}, e_{\left(1,2,2\right)} \pmod{{}_\chi M}.
\]
Suppose that for some $k_1,\ldots,k_4\in K$ we have
\begin{equation}
k_1e_{\left(1,1,1\right)}+k_2e_{\left(2,2,2\right)}+
k_3e_{\left(1,1,2\right)}+k_4e_{\left(1,2,2\right)}\in {}_\chi
M.\label{II.12}
\end{equation}

Applying the operator of $\chi$-symmetry $A_\chi=\sum_{\sigma\in
W}\chi^2(\sigma)\sigma$, we obtain
\[
k_1A_\chi e_{\left(1,1,1\right)}+k_2A_\chi e_{\left(2,2,2\right)}+
k_3A_\chi e_{\left(1,1,2\right)}+k_4A_\chi
e_{\left(1,2,2\right)}=0.
\]
On the other hand, $A_\chi e_{\left(1,1,1\right)}=A_\chi
e_{\left(2,2,2\right)}=0$, and  $A_\chi e_{\left(1,1,2\right)}$
and $A_\chi e_{\left(1,2,2\right)}$ are linearly independent over
$K$, hence $k_3=k_4=0$. Thus,
\[
k_1e_{\left(1,1,1\right)}+k_2e_{\left(2,2,2\right)}=
\ell_1(1-\varepsilon)e_{\left(1,1,1\right)}+\ell_2(1-\varepsilon)
e_{\left(2,2,2\right)} + f,
\]
where $\ell_1,\ell_2\in K$, and $f$ is a $K$-linear combination of
the tensors $e_{\left(1,1,2\right)}-\varepsilon
e_{\left(2,1,1\right)}$, $e_{\left(1,1,2\right)}-\varepsilon^2
e_{\left(1,2,1\right)}$, $e_{\left(1,2,2\right)}-\varepsilon
e_{\left(2,1,2\right)}$, and $e_{\left(1,2,2\right)}-\varepsilon^2
e_{\left(2,2,1\right)}$, that is, $k_1\in(1-\varepsilon)K$,
$k_2\in(1-\varepsilon)K$, and $f=0$. Therefore,~(\ref{II.12}) is
equivalent to $k_1\in(1-\varepsilon)K$, $k_2\in(1-\varepsilon)K$,
and $k_3=k_4=0$. In particular, the $K$-module $[\chi^2]^3(E)$ has
non-zero torsion part, hence it is not free.

\end{example}

 \section{Duality}

 \label{III}

Let the ring $K$ be an integral domain. Let us denote by $\cal F$
the category of $K$-modules with finite bases and, as usual,
denote by $Ob(\cal F)$ its set of objects. Let $E$ be a $K$-module
with finite basis $(e_\ell)_{\ell=1}^n$ and let $E^*$ be the dual
$K$-module with dual basis $(e_\ell^*)_{\ell=1}^n$. Denote by
$\langle\ ,\ \rangle$ the canonical bilinear form $E\times E^*\to
K$, $(x, x^*)\mapsto x^*(x)$. Let $W\leq S_d$ be a permutation
group with $|W|\in U(K)$, and let $\chi$ be a linear $K$-valued
character of $W$. We set $|W_\emptyset|=1$. For any $d\times
d$-matrix $A=(a_{ij})$ over $K$, the expression
\[
d_\chi^W(A)=\sum_{\sigma\in W}\chi(\sigma)a_{{\sigma^{-1}\left(1\right)}1}\ldots
a_{{\sigma^{-1}\left(d\right)}d}
\]
is known as (generalized) Schur function. It was introduced by
I.~Schur in \cite{[6]}.

\begin{theorem} \label{III.1} (i) The formulae
\begin{equation}
[\chi]^d(E)\times [\chi^{-1}]^d(E^*)\to  K, \label{III.2}
\end{equation}
\[
B(x_1\chi\ldots\chi x_d, x_1^*\chi^{-1}\ldots\chi^{-1} x_d^*)=d_\chi^W((<x_i,x_j^*>)_{i,j=1}^d),
\]
for $d\geq 1$, and the formula
\begin{equation}
[\chi]^0(E)\times [\chi^{-1}]^0(E^*)\to  K, \label{III.3}
\end{equation}
\[
B(k, k^*)=kk^*,
\]
define non-singular bilinear forms;

(ii) if $\iota_E^{\left(d\right)}\colon [\chi^{-1}]^d(E^*)\to
([\chi]^d(E))^*$ (resp., $\iota^{\left(0\right)}\colon
[\chi^{-1}]^0(E^*)\to ([\chi]^0(E))^*$) is the isomorphism of
$K$-modules, associated with~(\ref{III.2}) (resp.,
with~(\ref{III.3})), then the family
$\iota^{\left(d\right)}=(\iota_E^{\left(d\right)})_{E\in
Ob\left(\cal F\right)}$ (resp., $\iota^{\left(0\right)}$) is an
isomorphism of functors, $\iota^{\left(d\right)}\colon
[\chi^{-1}]^d(-^*)\to ([\chi]^d(-))^*$;

(iii) after the identifications via the functor
$\iota^{\left(d\right)}$ from (ii), $B$ is the canonical bilinear
form of the $K$-module $[\chi]^d(E)$, and the bases $(e_j)_{j\in
J}$ and $((1/|W_j|)e_j^*)_{j\in J}$ are dual.

\end{theorem}

\begin{proof} (i) For $d=0$ we get the multiplication of the ring
$K$. Let us suppose $d\geq 1$. The product $E^d\times (E^*)^d$
has a natural structure of $W\times W$-module (see \cite[2.1]{[3]}),
and the map
\[
E^d\times (E^*)^d\to K,
\]
\[
(x_1,\ldots,x_d, x_1^*,\ldots, x_d^*)
\mapsto  d_\chi^W((<x_i,x_j^*>)_{i,j=1}^d),
\]
is semi-symmetric of weight $\chi$ with respect to variables
$x_1,\ldots ,x_d$, and semi-symmetric of weight $\chi^{-1}$ with
respect to variables $x_1^*,\ldots ,x_d^*$. Hence by \cite[Lemma
2.1.2]{[3]} it gives rise to a bilinear form $B$ given by
formulae~(\ref{III.2}). We have $J(\chi,n,p)=J(\chi^{-1},n,p)=J$,
and in accord with Corollary~\ref{II.7}, $(e_j)_{j\in J}$ is a
basis for $[\chi]^d(E)$, and $(e_j^*)_{j\in J}$ is a basis for
$[\chi^{-1}]^d(E^*)$. If $\delta(j,k)$ is Kronecker's delta, then
\[
B(e_j,e_k^*)=\sum_{\sigma\in W}\chi^{-1}(\sigma)
\langle e_{j_{\sigma^{-1}\left(1\right)}},e_{k_1}^*\rangle\ldots
\langle e_{j_{\sigma^{-1}\left(d\right)}},e_{k_d}^*\rangle
\]
\[
=\sum_{\sigma\in W}\chi^{-1}(\sigma)
\delta(j_{\sigma^{-1}\left(1\right)},k_1)\ldots
\delta(j_{\sigma^{-1}\left(d\right)},k_d),
\]
hence
\begin{equation}
B(e_j,e_k^*)=|W_j|\delta(j,k).\label{III.4}
\end{equation}

In particular,~(\ref{III.2}) and~(\ref{III.3}) are non-singular
forms for any $d\geq 0$.

(ii) For any $K$-linear map $u\colon E\to F$ we denote by
${}^t\hskip-2ptu\colon F^*\to E^*$ its transpose. A direct
computation shows that
\begin{equation}
{}^t\hskip-2pt([\chi]^d(u))\circ
\iota_F^{\left(d\right)}=\iota_E^{\left(d\right)}\circ
([\chi^{-1}]^d({}^t\hskip-2pt u)).\label{III.5}
\end{equation}

(iii) The equality~(\ref{III.4}) yields that $(e_j)_{j\in J}$,
$(\frac{1}{|W_j|}e_j^*)_{j\in J}$ is a pair of dual bases.

\end{proof}

\begin{remark} \label{III.6} Throughout the end of the paper we will use notation
\[
\langle x_1\chi\ldots\chi x_d, x_1^*\chi^{-1}\ldots\chi^{-1} x_d^*\rangle=
B(x_1\chi\ldots\chi x_d, x_1^*\chi^{-1}\ldots\chi^{-1} x_d^*),
\]
and in this notation, for any $x=\sum_{j\in J}x_je_j$, and for any
$x^*=\sum_{j\in J}x_j^*e_j^*$, one has
\begin{equation}
\langle x, x^*\rangle=\sum_{j\in J}|W_j|x_jx_j^*.\label{III.7}
\end{equation}

\end{remark}

\begin{remark} \label{III.8} In accord with Theorem~\ref{III.1}, (ii), (iii),
for any $d\geq 1$, and for any $K$-module $E$ with finite basis we
identify $([\chi]^d(E))^*$ with $[\chi^{-1}]^d(E^*)$ as
$K$-modules via the functor $\iota^{\left(d\right)}$, and call the
elements of $[\chi^{-1}]^d(E^*)$ $d-\chi$-forms on $E$.

\end{remark}

\begin{corollary} \label{III.9} For any $K$-linear map $u\colon E\to F$ one has
${}^t\hskip-2pt([\chi]^du)=[\chi^{-1}]^d({}^t\hskip-2pt u)$.

\end{corollary}

\begin{proof} This the equality~(\ref{III.5}) after
the identifications via the functor $\iota^{\left(d\right)}$.

\end{proof}

Let $A=(a_{r,s})$ be an $m\times n$ matrix  over $K$ and let
$d\geq 1$. For any $j\in J(\chi,m,d)$, $k\in J(\chi,n,d)$, we set
$a_{jk}=\prod_{t=1}^da_{j_tk_t}$, and
\[
A_{\left(j\right)k}(\chi)=\sum_{\tau\in W^{\left(j\right)}}\chi^{-1}(\tau)a_{\tau jk},
\]
and call the expression $A_{\left(j\right)k}(\chi)$ the
\emph{$(j,k)$-th row minor of weight $\chi$} of $A$.

Let $A=(a_{rt})$ and $A^\prime=(a_{sh}^\prime)$ be two $n\times d$
matrices over $K$. Using notation from the beginning of
Section~\ref{III}, we set $x_t=\sum_{r=1}^na_{rt}e_r$,
$x_h^*=\sum_{s=1}^na_{sh}^\prime e_s^*$, where $t,h=1,\ldots ,d$.
Then $\langle x_t, x_h^*\rangle=\sum_{r=1}^na_{rt}a_{rh}^\prime$
is the $th$-entry of the matrix ${}^t\hskip-2pt AA^\prime$, and
hence
\begin{equation}
\langle x_1\chi\ldots\chi x_d, x_1^*\chi^{-1}\ldots\chi^{-1}
x_d^*\rangle=d_\chi({}^t\hskip-2ptAA^\prime).\label{III.10}
\end{equation}

On the other hand,
\begin{equation}
x_1\chi\ldots\chi x_d=\sum_{j\in
J}A_{\left(j\right)}(\chi)e_j,\mbox{\ }
x_1^*\chi^{-1}\ldots\chi^{-1} x_d^*=\sum_{j\in
J}A_{\left(j\right)}^\prime(\chi^{-1})e_j^*,\label{III.11}
\end{equation}
where $A_{\left(j\right)}(\chi)=A_{\left(j\right)k}(\chi)$, and
$A_{\left(j\right)}^\prime(\chi^{-1})=A_{\left(j\right)k}^\prime(\chi^{-1})$
with $k=(1,...,d)$. Therefore~(\ref{III.7}) and~(\ref{III.10})
yield
\[
d_\chi({}^t\hskip-2pt AA^\prime)=\sum_{j\in
J}|W_j|A_{\left(j\right)}(\chi)A_{\left(j\right)}^\prime(\chi^{-1}).
\]
In particular, when $A=A'$ we obtain generalized Lagrange identity
\[
d_\chi({}^t\hskip-2pt AA)=\sum_{j\in J}|W_j|A_{\left(j\right)}(\chi)A_{\left(j\right)}(\chi^{-1}).
\]

\begin{lemma} \label{III.12} Let $A=(a_{th})$ be a $d\times d$ matrix over $K$. Then, in
the previous notations, one has:

(i) $d_\chi({}^tA)=d_{\chi^{-1}}(A)$;

(ii) $d_\chi(A)=\langle x_1\chi^{-1}\ldots\chi^{-1}x_d,
e_1^*\chi\ldots\chi e_d^*\rangle$;

(iii) The generalized Schur function $d_\chi(A)$ is semi-symmetric
of weight $\chi^{-1}$ (resp., of weight $\chi$) with respect to
the columns (resp., the rows) of the matrix $A$.

\end{lemma}

\begin{proof} (i) Direct checking.

(ii)  Using (i) and~(\ref{III.10}) with $d=n$ and $A^\prime=I_d$
(the unit $d\times d$ matrix), we obtain the equality.

(iii) This is an immediate consequence of (ii) and (i).

\end{proof}

Throughout the end of the paper we fix the following notation:  

$K$ is both a $\hbox{\ccc Q}$-ring and an integral domain;

$(\chi_d\colon W_d\to K)_{d\geq 1}$ is an $\omega$-invariant
sequence of characters; 

$E$ is a $K$-module with finite basis;

$[\chi](E)$ is the semi-symmetric algebra of weight $\chi$ of $E$.  

We remind that the dual graded $K$-module
$([\chi](E))^{*gr}$ is, by definition, the direct sum
$\coprod_{d\geq 0}([\chi]^d(E))^*$, where we identify a linear
form on $[\chi]^d(E)$ with its extension by $0$ to $[\chi](E)$.
Let us set $\iota=\coprod_{d\geq 0}\iota^{\left(d\right)}$.

 Since the $K$-module $E$ has a finite basis, then it is a projective
 module of finite type, and using Corollary~\ref{II.8}, 
 \cite[A II, p. 80, Cor. 1]{[2]}, and
 Theorem~\ref{III.1}, we obtain

\begin{theorem} \label{III.13}(i) $\iota\colon [\chi](-^*)\to ([\chi](-))^{*gr}$
is an isomorphism of functors;

(ii) After the identification via the functor $\iota$ from (i),
the restriction of the canonical bilinear form of the $K$-module
$[\chi](E)$ on $[\chi](E)\times [\chi](E^*)$ is given by the
formulae
\begin{equation}
\langle\ ,\ \rangle\colon [\chi](E)\times [\chi](E^*)\to
K,\label{III.14}
\end{equation}
\[
 \langle x_1\chi\ldots\chi x_r, x_1^*\chi\ldots\chi x_s^*\rangle=
  \left\{
\begin{array}{ll}
0 & \mbox{if $r\neq s$}\\
d_\chi((\langle x_i,x_j^*\rangle)_{i,j=1}^r) & \mbox{if $r=s\geq 1$}\\
1 & \mbox{if $r=s=0$};
\end{array}
\right.
\]

(iii) for any $k\geq 2$ the restriction of the canonical bilinear
form of the $K$-module $([\chi](E))^{\otimes k}$ on
$([\chi](E))^{\otimes k}\times [\chi](E^*))^{\otimes k}$ is given
by the formulae
\begin{equation}
\langle\ ,\ \rangle\colon ([\chi](E))^{\otimes k}\times
[\chi](E^*))^{\otimes k}\to K,\label{III.15}
\end{equation}
\[
 \langle x_1\chi\ldots\chi x_r\otimes x_1\chi\ldots\chi x_{r'}\otimes\cdots,
 x_1^*\chi\ldots\chi x_s^*\otimes x_1^*\chi\ldots\chi x_{s'}^*\otimes\cdots\rangle=
 \]
 \[
  \left\{
\begin{array}{ll}
0 & \mbox{if $(r,r',\ldots)\neq (s,s',\ldots)$}\\
\langle x_1\chi\ldots\chi x_r, x_1^*\chi\ldots\chi x_r^*\rangle
\langle x_1\chi\ldots\chi x_{r'}, x_1^*\chi\ldots\chi
x_{r'}^*\rangle\ldots & \mbox{if $(r,r',\ldots)=(s,s',\ldots)$},
\end{array}
\right.
\]

\end{theorem}

\begin{remark}\label{III.16} Let $d$, $e\ldots$, $h$, be non-negative integers with
$d+e+\cdots+h=n$. We set
\[
J(\chi;n;d,e,\ldots, h)=
\]
\[
\{(j,k,\ldots, r)\in J(\chi,n,d)\times
J(\chi,n,e)\times\cdots\times J(\chi,n,h)\mid \ell m(j,k,\ldots,r)=(1,...,n)\}.
\]
Let $M(\chi;n;d,e,\ldots, h)$ be the set of lexicographically
minimal representatives of left classes of $W_n$ modulo
$W_d\times\omega^d(W_e)\times\cdots\times\omega^{d+e+\cdots}(W_h)$.
We identify the set $M(\chi;n;d,e,\ldots, h)$ with the set
$J(\chi;n;d,e,\ldots, h)$ via the canonical bijection
\[
M(\chi;n;d,e,\ldots, h)\to J(\chi;n;d,e,\ldots, h),
\]
\[
\zeta\mapsto
((\zeta(1),\ldots,\zeta(d)),(\zeta(d+1),\ldots,\zeta(d+e)),
\ldots, (\zeta(d+e+\cdots+1),\ldots,\zeta(n))).
\]

\end{remark}

We fix $(\lambda,\mu,\ldots, \nu)\in J(\chi;n;d,e,\ldots, h)$, and
let $\sigma\in W_n$ be a permutation, such that
$\lambda=(\sigma(1),\ldots ,\sigma(d))$, $\mu=(\sigma(d+1),\ldots
,\sigma(d+e))$, $\ldots$, and $\nu=(\sigma(d+e+\cdots+1),\ldots
,\sigma(n))$. We have $\zeta(\lambda,\mu,\ldots,
\nu)=\chi(\sigma)$. Let us write $d_\chi(A)$ for $d_{\chi_n}(A)$.

\begin{proposition} \label{III.17} Let $A$ be an $n\times n$ matrix over $K$. Then
\[
d_\chi(A)=
\]
\[
\zeta(\lambda,\mu,\ldots, \nu)\sum_{\left(j,k,\ldots, r\right)\in J\left(\chi;n;d,e,\ldots, h\right)}
\zeta(j,k,\ldots, r)A_{\left(j\right)\lambda}(\chi)A_{\left(k\right)\mu}(\chi)\ldots A_{\left(r\right)\nu}(\chi).
\]
(Laplace expansion of $d_\chi(A)$ with respect to $\lambda$, $\mu$, $\ldots$, $\nu$).

\end{proposition}

\begin{proof} Indeed, using Lemma~\ref{III.12}, (ii), Corollary~\ref{I.3}, (i),
the expansions~(\ref{III.11}), and Corollary~\ref{II.10}, (ii), we
obtain
\[
d_\chi(A)=\langle x_1\chi^{-1}\ldots\chi^{-1}x_n,
e_1^*\chi\ldots\chi e_n^*\rangle= \langle x_1\chi\ldots\chi x_n,
e_1^*\chi\ldots\chi e_n^*\rangle=
\]
\[
\zeta(\lambda,\mu,\ldots, \nu)\langle x_{\lambda_1}\chi\ldots\chi x_{\lambda_d}\chi x_{\mu_1}\chi
\ldots\chi x_{\mu_e}\chi x_{\nu_1}\chi\ldots\chi x_{\nu_h}, e_1^*\chi\ldots\chi e_n^*\rangle=
\]
\[
\zeta(\lambda,\mu,\ldots, \nu)\langle
\sum_{\left(j,k,\ldots, r\right)\in J\left(\chi,n,d\right)\times
J\left(\chi,n,e\right)\times\cdots\times J\left(\chi,;n,h\right)}
\zeta(j,k,\ldots, r)
\]
\[
A_{\left(j\right)\lambda}(\chi)A_{\left(k\right)\mu}(\chi)\ldots
A_{\left(r\right)\nu}(\chi)e_{\ell m\left(j,k,\ldots, r\right)},
e_1^*\chi\ldots\chi e_n^*\rangle=
 \]
 \[
 \zeta(\lambda,\mu,\ldots\nu)\sum_{\left(j,k,\ldots, r\right)\in J\left(\chi;n;d,e,\ldots,h\right)}
\zeta(j,k,\ldots,
r)A_{\left(j\right)\lambda}(\chi)A_{\left(k\right)\mu}(\chi)\ldots
A_{\left(r\right)\nu}(\chi).
\]

\end{proof}

\begin{proposition}\label{III.18} For any non-negative integers
$d$, $e$,$\ldots$, $h$ with $d+e+\cdots +h=n$ one has the
following expansions of the bilinear form~(\ref{III.14}) (Laplace
expansions):
\[
\langle x_1\chi\ldots\chi x_n, x_1^*\chi\ldots\chi x_n^*\rangle=
\]
\[
\sum_{\zeta\in M\left(\chi;n;d,e,\ldots, h\right)}
\chi(\zeta)
\langle x_{\zeta\left(1\right)}\chi\ldots\chi x_{\zeta\left(d\right)}, x_1^*\chi
\ldots\chi x_d^*\rangle
\]
\[
\langle x_{\zeta\left(d+1\right)}\chi\ldots\chi x_{\zeta\left(d+e\right)}, x_{d+1}^*\chi
\ldots\chi x_{d+e}^*\rangle
\]
\[
\cdots\langle x_{\zeta\left(d+e+\cdots+1\right)}\chi\ldots\chi
x_{\zeta\left(n\right)}, x_{d+e+\cdots+1}^*\chi \ldots\chi
x_n^*\rangle=
\]
\[
\sum_{\zeta\in M\left(\chi;n;d,e,\ldots, h\right)} \chi(\zeta)
\langle x_1\chi \ldots\chi x_d,
x_{\zeta\left(1\right)}^*\chi\ldots\chi
x_{\zeta\left(d\right)}^*\rangle
\]
\[
\langle x_{d+1}\chi \ldots\chi
x_{d+e},x_{\zeta\left(d+1\right)}^*\chi\ldots\chi
x_{\zeta\left(d+e\right)}^*\rangle
\]
\[
\cdots\langle x_{d+e+\cdots+1}\chi \ldots\chi x_n
,x_{\zeta\left(d+e+\cdots+1\right)}^*\chi\ldots\chi
x_{\zeta\left(n\right)}^*\rangle.
\]

\end{proposition}

\begin{proof} We have
\[
\langle x_1\chi\ldots\chi x_n, x_1^*\chi\ldots\chi x_n^*\rangle=
\]
\[
\sum_{\zeta'\in W_n}\chi(\zeta')
(\langle x_{\zeta'\left(1\right)}, x_1^*\rangle
\cdots \langle x_{\zeta'\left(d\right)}, x_d^*\rangle)
(\langle x_{\zeta'\left(d+1\right)}, x_{d+1}^*\rangle
\cdots \langle x_{\zeta'\left(d+e\right)}, x_{d+e}^*\rangle)
\]
\[
\cdots (\langle x_{\zeta'\left(d+e+\cdots+1\right)}, x_{d+e+\cdots+1}^*\rangle
\cdots \langle x_{\zeta'\left(n\right)}, x_n^*\rangle)=
\]
\[
\sum_{\zeta\in M\left(\chi;n;d,e,\ldots, h\right)}
\sum_{\left(\sigma',\tau',\ldots,\eta'\right)\in W_d\times\omega^d\left(W_e\right)
\times\cdots\times\omega^{d+e+\cdots}\left(W_h\right)}
\chi(\zeta)\chi(\sigma')\chi(\tau')\ldots \chi(\eta')
\]
\[
(\langle x_{\zeta\left(\sigma'\left(1\right)\right)}, x_1^*\rangle
\cdots \langle x_{\zeta\left(\sigma'\left(d\right)\right)}, x_d^*\rangle)
(\langle x_{\zeta\left(\tau'\left(d+1\right)\right)}, x_{d+1}^*\rangle
\cdots \langle x_{\zeta\left(\tau'\left(d+e\right)\right)}, x_{d+e}^*\rangle)
\]
\[
\cdots (\langle x_{\zeta\left(\eta'\left(d+e+\cdots+1\right)\right)},
x_{d+e+\cdots+1}^*\rangle
\cdots \langle x_{\zeta\left(\eta'\left(n\right)\right)}, x_n^*\rangle)=
\]
\[
\sum_{\zeta\in M\left(\chi;n;d,e,\ldots, h\right)}
\chi(\zeta)
(\sum_{\sigma'\in W_d}\chi(\sigma')
\langle x_{\zeta\left(\sigma'\left(1\right)\right)}, x_1^*\rangle
\cdots \langle x_{\zeta\left(\sigma'\left(d\right)\right)}, x_d^*\rangle)
\]
\[
(\sum_{\tau'\in \omega^d\left(W_e\right)}\chi(\tau')
\langle x_{\zeta\left(\tau'\left(d+1\right)\right)}, x_{d+1}^*\rangle
\cdots \langle x_{\zeta\left(\tau'\left(d+e\right)\right)}, x_{d+e}^*\rangle)
\]
\[
\cdots (\sum_{\eta'\in \omega^{d+e+\cdots}\left(W_d\right)}\chi(\eta')
\langle x_{\zeta\left(\eta'\left(d+e+\cdots+1\right)\right)}, x_{d+e+\cdots+1}^*\rangle
\cdots \langle x_{\zeta\left(\eta'\left(n\right)\right)}, x_n^*\rangle)=
\]
\[
\sum_{\zeta\in M\left(\chi;n;d,e,\ldots, h\right)}
\chi(\zeta)
\langle x_{\zeta\left(1\right)}\chi\ldots\chi x_{\zeta\left(d\right)}, x_1^*\chi
\ldots\chi x_d^*\rangle
\]
\[
\langle x_{\zeta\left(d+1\right)}\chi\ldots\chi
x_{\zeta\left(d+e\right)}, x_{d+1}^*\chi \ldots\chi
x_{d+e}^*\rangle\cdots
\]
\[
\langle x_{\zeta\left(d+e+\cdots+1\right)}\chi\ldots\chi
x_{\zeta\left(n\right)}, x_{d+e+\cdots+1}^*\chi\ldots\chi
x_n^*\rangle.
\]
For the second equality, we can write
\[
\langle x_1\chi\ldots\chi x_n, x_1^*\chi\ldots\chi x_n^*\rangle=
\]
\[
\sum_{\zeta'\in W_n}\chi(\zeta') (\langle x_1,
x_{\zeta'^{-1}\left(1\right)}^*\rangle \cdots \langle x_d,
x_{\zeta'^{-1}\left(d\right)}^*\rangle) (\langle x_{d+1}
,x_{\zeta'^{-1}\left(d+1\right)}^*\rangle \cdots \langle x_{d+e},
x_{\zeta'^{-1}\left(d+e\right)}^*\rangle)
\]
\[
\cdots (\langle x_{d+e+\cdots+1},
x_{\zeta'^{-1}\left(d+e+\cdots+1\right)}^*\rangle \cdots \langle
x_n, x_{\zeta'^{-1}\left(n\right)}^*\rangle)=
\]
\[
\sum_{\zeta'\in W_n}\chi(\zeta') (\langle x_1,
x_{\zeta'\left(1\right)}^*\rangle \cdots \langle x_d,
x_{\zeta'\left(d\right)}^*\rangle) (\langle x_{d+1}
,x_{\zeta'\left(d+1\right)}^*\rangle \cdots \langle x_{d+e},
x_{\zeta'\left(d+e\right)}^*\rangle)
\]
\[
\cdots (\langle x_{d+e+\cdots+1},
x_{\zeta'\left(d+e+\cdots+1\right)}^*\rangle \cdots \langle x_n,
x_{\zeta'\left(n\right)}^*\rangle),
\]
and then we proceed by analogy.

\end{proof}

According to Lemma~\ref{A.4}, for any $n$ we obtain a $K$-linear map
\[
[\chi]^n(E)\to\oplus_{d+e+\cdots+h=n}[\chi]^d(E)\hskip-2pt\otimes\hskip-2pt
[\chi]^e(E)\ldots\otimes\hskip-2pt [\chi]^h(E),
\]
\[
x_1\chi\ldots\chi x_n\mapsto \sum_{d+e+\cdots+h=n}\sum_{\rho\in
M\left(\chi;n;d,e,\ldots, h\right)}\chi(\rho)
\]
\[ (x_{\rho\left(1\right)}\chi\ldots\chi
x_{\rho\left(d\right)})\hskip-2pt\otimes\hskip-2pt
(x_{\rho\left(d+1\right)}\chi\ldots\chi
x_{\rho\left(d+e\right)})\otimes\ldots\otimes\hskip-2pt
(x_{\rho\left(d+e+\cdots+1\right)}\chi\ldots\chi
x_{\rho\left(n\right)}).
\]
Therefore, for any $k\geq 2$ we get a homomorphism of graded
$K$-modules
\begin{equation}
c_k(E)\colon [\chi](E)\to ([\chi](E))^{\otimes k},\label{III.19}
\end{equation}
\[
c_k(E)(x_1\chi\ldots\chi x_n)= \sum_{d+e+\cdots+h=n}\sum_{\rho\in
M\left(\chi;n;d,e,\ldots, h\right)}\chi(\rho)
\]
\[(x_{\rho\left(1\right)}\chi\ldots\chi
x_{\rho\left(d\right)})\hskip-2pt\otimes\hskip-2pt
(x_{\rho\left(d+1\right)}\chi\ldots\chi
x_{\rho\left(d+e\right)})\hskip-2pt\otimes\ldots\otimes\hskip-2pt
(x_{\rho\left(d+e+\cdots+1\right)}\chi\ldots\chi
x_{\rho\left(n\right)}).
\]

\begin{corollary}\label{III.20} For any $k$ in number non-negative integers
$d$, $e$,$\ldots$, $h$ with $d+e+\cdots +h=n$ one has
\[
\langle x_1\chi\ldots\chi x_n, x_1^*\chi\ldots\chi x_n^*\rangle=
\]
\[
\langle c_k(E)(x_1\chi\ldots\chi x_n), x_1^*\chi \ldots\chi
x_d^*\hskip-2pt\otimes\hskip-2pt x_{d+1}^*\chi \ldots\chi
x_{d+e}^*\hskip-2pt\otimes\cdots\otimes\hskip-2pt
x_{d+e+\cdots+1}^*\chi \ldots\chi x_n^*\rangle=
\]
\[
\langle x_1\chi \ldots\chi x_d\hskip-2pt\otimes\hskip-2pt
x_{d+1}\chi \ldots\chi x_{d+e}\hskip-2pt\otimes\cdots\otimes
\hskip-2ptx_{d+e+\cdots+1}\chi \ldots\chi x_n, c_k(E^*)(x_1^*\chi
\ldots\chi x_n^*)\rangle.
\]

\end{corollary}

\begin{proof} Using~(\ref{III.15}), and Proposition~\ref{III.18}, we have
\[
\langle x_1\chi\ldots\chi x_n, x_1^*\chi\ldots\chi x_n^*\rangle=
\]
\[
\sum_{\zeta\in M\left(\chi;n;d,e,\ldots, h\right)} \chi(\zeta)
\langle x_{\zeta\left(1\right)}\chi\ldots\chi
x_{\zeta\left(d\right)}, x_1^*\chi \ldots\chi x_d^*\rangle
\]
\[
\langle x_{\zeta\left(d+1\right)}\chi\ldots\chi
x_{\zeta\left(d+e\right)}, x_{d+1}^*\chi \ldots\chi
x_{d+e}^*\rangle
\]
\[
\cdots\langle x_{\zeta\left(d+e+\cdots+1\right)}\chi\ldots\chi
x_{\zeta\left(n\right)}, x_{d+e+\cdots+1}^*\chi \ldots\chi
x_n^*\rangle=
\]
\[
\sum_{\zeta\in M\left(\chi;n;d,e,\ldots, h\right)} \chi(\zeta)
\]
\[
\langle x_{\zeta\left(1\right)}\chi\ldots\chi
x_{\zeta\left(d\right)}\hskip-2pt\otimes \hskip-2pt
x_{\zeta\left(d+1\right)}\chi\ldots\chi
x_{\zeta\left(d+e\right)}\hskip-2pt\otimes\cdots\otimes \hskip-2pt
x_{\zeta\left(d+e+\cdots+1\right)}\chi\ldots\chi
x_{\zeta\left(n\right)},
\]
\[
x_1^*\chi \ldots\chi x_d^*\hskip-2pt\otimes\hskip-2pt
x_{d+1}^*\chi \ldots\chi x_{d+e}^*\hskip-2pt\otimes\cdots\otimes
\hskip-2pt x_{d+e+\cdots+1}^*\chi \ldots\chi x_n^*\rangle =
\]
\[
\langle c_k(E)(x_1\chi\ldots\chi x_n), x_1^*\chi \ldots\chi
x_d^*\hskip-2pt\otimes\hskip-2pt x_{d+1}^*\chi \ldots\chi
x_{d+e}^*\hskip-2pt\otimes\cdots\otimes\hskip-2pt
x_{d+e+\cdots+1}^*\chi \ldots\chi x_n^*\rangle.
\]

Similarly, using the second identity of Proposition~\ref{III.18},
we obtain the second identity of this corollary.

\end{proof}

\section{Coalgebra properties}

\label{IV}

Let us set $c_k=c_k(E)$, and $c_E=c_2(E)$, where $c_k(E)$, $k\geq
2$, is the homomorphism of graded $K$-modules from~(\ref{III.19}).

\begin{proposition} \label{IV.1} One has
\[
c_k=(c_{k-1}\hskip-2pt\otimes\hskip-2pt 1)\circ
c_E=(1\hskip-2pt\otimes\hskip-2pt c_{k-1})\circ c_E,
\]
where $1$ is the identity map of $[\chi](E)$. 

\end{proposition}

\begin{proof} We have
\[
c_E(x_1\chi\ldots\chi x_n)= \sum_{p+h=n}\sum_{\rho\in
M\left(\chi;n;p,h\right)}\chi(\rho)
\]
\[
(x_{\rho\left(1\right)}\chi\ldots\chi
x_{\rho\left(p\right)})\hskip-2pt\otimes\hskip-2pt
(x_{\rho\left(p+1\right)}\chi\ldots\chi x_{\rho\left(n\right)}).
\]
First, we apply the $K$-linear map
$c_{k-1}\hskip-2pt\otimes\hskip-2pt 1$ and get
\[
(c_{k-1}\hskip-2pt\otimes\hskip-2pt 1)(c_E(x_1\chi\ldots\chi
x_n))= \sum_{p+h=n}\sum_{\rho\in
M\left(\chi;n;p,h\right)}\chi(\rho)
\]
\[
c_{k-1}(x_{\rho\left(1\right)}\chi\ldots\chi
x_{\rho\left(p\right)})\hskip-2pt\otimes\hskip-2pt
(x_{\rho\left(p+1\right)}\chi\ldots\chi x_{\rho\left(n\right)})=
\]
\[
\sum_{p+h=n}\sum_{\rho\in M\left(\chi;n;p,h\right)}
\sum_{d+e+\cdots=p}\sum_{\varrho\in
M\left(\chi;p;d,e,\ldots\right)}\chi(\rho\varrho)
(x_{\rho\left(\varrho\left(1\right)\right)}\chi\ldots\chi
x_{\rho\left(\varrho\left(d\right)\right)})
\]
\[
\hskip-2pt\otimes\hskip-1pt
(x_{\rho\left(\varrho\left(d+1\right)\right)}\chi\ldots\chi
x_{\rho\left(\varrho\left(d+e\right)\right)})\hskip-2pt
\otimes\ldots\otimes\hskip-2pt
(x_{\rho\left(p+1\right)}\chi\ldots\chi x_{\rho\left(n\right)})=
\]
\[
\sum_{p+h=n}\sum_{\rho\in M\left(\chi;n;p,h\right)}
\sum_{d+e+\cdots=p}\sum_{\varrho\in
M\left(\chi;p;d,e,\ldots\right)}\chi(\rho\varrho)
(x_{\rho\left(\varrho\left(1\right)\right)}\chi\ldots\chi
x_{\rho\left(\varrho\left(d\right)\right)})
\]
\[
\hskip-2pt\otimes\hskip-1pt
(x_{\rho\left(\varrho\left(d+1\right)\right)}\chi\ldots\chi
x_{\rho\left(\varrho\left(d+e\right)\right)})\hskip-2pt
\otimes\ldots\otimes\hskip-2pt
(x_{\rho\left(\varrho\left(p+1\right)\right)}\chi\ldots\chi
x_{\rho\left(\varrho\left(n\right)\right)})=
\]
\[
\sum_{d+e+\cdots+h=n}\sum_{\left(\rho,\varrho\right)\in
M\left(\chi;n;p,h\right)\times M\left(\chi;p;d,e,\ldots\right)}
\chi(\rho\varrho)
(x_{\rho\left(\varrho\left(1\right)\right)}\chi\ldots\chi
x_{\rho\left(\varrho\left(d\right)\right)})
\]
\[
\hskip-1pt\otimes\hskip-1pt
(x_{\rho\left(\varrho\left(d+1\right)\right)}\chi\ldots\chi
x_{\rho\left(\varrho\left(d+e\right)\right)})\hskip-2pt
\otimes\ldots\otimes\hskip-2pt
(x_{\rho\left(\varrho\left(p+1\right)\right)}\chi\ldots\chi
x_{\rho\left(\varrho\left(n\right)\right)}).
\]
In terms of Notation~\ref{A.2}, we set
$\rho\varrho\sigma'\tau'\ldots\eta'=1\cdot(\rho\varrho)$, where
$\sigma'\in W_d$, $\tau'\in \omega^d(W_e)$,$\ldots$, $\eta'\in
\omega^p(W_h)$, $\sigma'=\sigma$, $\sigma\in W_d$,
$\tau'=\omega^d(\tau)$, $\tau\in W_e$,$\ldots$,
$\eta'=\omega^p(\eta)$, $\eta\in W_h$. Then
$\chi(\rho\varrho)\chi(\sigma)\chi(\tau)\ldots\chi(\eta)=\chi(1\cdot(\rho\varrho))$,
and we have
\[
(c_{k-1}\hskip-2pt\otimes\hskip-2pt 1)(c_E(x_1\chi\ldots\chi
x_n))=
\]
\[
\sum_{d+e+\cdots+h=n}\sum_{\left(\rho,\varrho\right)\in
M\left(\chi;n;p,h\right)\times M\left(\chi;p;d,e,\ldots\right)}
\chi(1\cdot(\rho\varrho))
(x_{\rho\left(\varrho\left(\sigma\left(1\right)\right)\right)}\chi\ldots\chi
x_{\rho\left(\varrho\left(\sigma\left(d\right)\right)\right)})
\]
\[
\hskip-1pt\otimes\hskip-1pt
(x_{\rho\left(\varrho\left(\tau\left(d+1\right)\right)\right)}\chi\ldots\chi
x_{\rho\left(\varrho\left(\tau\left(d+e\right)\right)\right)})\hskip-2pt
\otimes\ldots\otimes\hskip-2pt
(x_{\rho\left(\varrho\left(\eta\left(p+1\right)\right)\right)}\chi\ldots\chi
x_{\rho\left(\varrho\left(n\right)\right)})=
\]
\[
\sum_{d+e+\cdots+h=n}\sum_{\left(\rho,\varrho\right)\in
M\left(\chi;n;p,h\right)\times M\left(\chi;p;d,e,\ldots\right)}
\chi(1\cdot(\rho\varrho))
(x_{\left(1\cdot\left(\rho\varrho\right)\right)
\left(1\right)}\chi\ldots\chi
x_{\left(1\cdot\left(\rho\varrho\right)\right)\left(d\right)})
\]
\[
\hskip-1pt\otimes\hskip-1pt
(x_{\left(1\cdot\left(\rho\varrho\right)\right)
\left(d+1\right)}\chi\ldots\chi
x_{\left(1\cdot\left(\rho\varrho\right)\right) \left(d+e\right)})
\hskip-2pt\otimes\ldots\otimes\hskip-2pt
(x_{\left(1\cdot\left(\rho\varrho\right)\right)
\left(p+1\right)}\chi\ldots\chi
x_{\left(1\cdot\left(\rho\varrho\right)\right) \left(n\right)}).
\]
According to Lemma~\ref{A.5} we obtain
\[
(c_{k-1}\hskip-2pt\otimes\hskip-2pt 1)(c_E(x_1\chi\ldots\chi
x_n))=
\]
\[
\sum_{d+e+\cdots+h=n}\sum_{\varsigma\in M\left(\chi;n;d,e,\ldots,
h \right)} \chi(\varsigma) (x_{\varsigma
\left(1\right)}\chi\ldots\chi x_{\varsigma\left(d\right)})
\]
\[
\hskip-1pt\otimes\hskip-1pt (x_{\varsigma
\left(d+1\right)}\chi\ldots\chi x_{\varsigma \left(d+e\right)})
\hskip-2pt\otimes\ldots\otimes\hskip-2pt (x_{\varsigma
\left(p+1\right)}\chi\ldots\chi x_{\varsigma\left(n\right)})=
\]
\[
c_k(x_1\chi\ldots\chi x_n).
\]
Similarly, we apply the $K$-linear map $1\hskip-2pt\otimes
\hskip-2pt c_{k-1}$ and obtain
\[
(1\hskip-2pt\otimes\hskip-2pt c_{k-1})(c_E(x_1\chi\ldots\chi
x_n))= \sum_{d+q=n}\sum_{\rho\in
M\left(\chi;n;d,q\right)}\chi(\rho)
\]
\[
(x_{\rho\left(1\right)}\chi\ldots\chi
x_{\rho\left(d\right)})\hskip-2pt\otimes\hskip-2pt
c_{k-1}(x_{\rho\left(d+1\right)}\chi\ldots\chi
x_{\rho\left(n\right)})=
\]
\[
\sum_{d+q=n}\sum_{\rho\in M\left(\chi;n;d,q\right)}
\sum_{e+\cdots+h=q}\sum_{\varrho\in
\omega^d\left(M\left(\chi;q;e,\ldots,h\right)\right)}\chi(\rho\varrho)
(x_{\rho\left(1\right)}\chi\ldots\chi x_{\rho\left(d\right)})
\]
\[
\hskip-2pt\otimes\hskip-1pt
(x_{\rho\left(\varrho\left(d+1\right)\right)}\chi\ldots\chi
x_{\rho\left(\varrho\left(d+e\right)\right)})\hskip-2pt
\otimes\ldots\otimes\hskip-2pt
(x_{\rho\left(\varrho\left(p+1\right) \right)}\chi\ldots\chi
x_{\rho\left(\varrho\left(n\right)\right)})=
\]
\[
\sum_{d+q=n}\sum_{\rho\in M\left(\chi;n;d,q\right)}
\sum_{e+\cdots+h=q}\sum_{\varrho\in
\omega^d\left(M\left(\chi;q;e,\ldots,h\right)\right)}\chi(\rho\varrho)
(x_{\rho\left(\varrho\left(1\right)\right)}\chi\ldots\chi
x_{\rho\left(\varrho\left(d\right)\right)})
\]
\[
\hskip-2pt\otimes\hskip-1pt
(x_{\rho\left(\varrho\left(d+1\right)\right)}\chi\ldots\chi
x_{\rho\left(\varrho\left(d+e\right)\right)})\hskip-2pt
\otimes\ldots\otimes\hskip-2pt
(x_{\rho\left(\varrho\left(p+1\right) \right)}\chi\ldots\chi
x_{\rho\left(\varrho\left(n\right)\right)})=
\]
\[
\sum_{d+e+\cdots+h=n}\sum_{\left(\rho,\varrho\right)\in
M\left(\chi;n;d,q\right)\times
\omega^d\left(M\left(\chi;q;e,\ldots,h\right)\right)}\chi(\rho\varrho)
(x_{\rho\left(\varrho\left(1\right)\right)}\chi\ldots\chi
x_{\rho\left(\varrho\left(d\right)\right)})
\]
\[
\hskip-2pt\otimes\hskip-1pt
(x_{\rho\left(\varrho\left(d+1\right)\right)}\chi\ldots\chi
x_{\rho\left(\varrho\left(d+e\right)\right)})
\otimes\ldots\otimes\hskip-2pt
(x_{\rho\left(\varrho\left(p+1\right) \right)}\chi\ldots\chi
x_{\rho\left(\varrho\left(n\right)\right)})=
\]
\[
\sum_{d+e+\cdots+h=n}\sum_{\left(\rho,\varrho\right)\in
M\left(\chi;n;d,q\right)\times
\omega^d\left(M\left(\chi;q;e,\ldots,h\right)\right)}
\chi(1\cdot(\rho\varrho))
(x_{\left(1\cdot\left(\rho\varrho\right)\right)
\left(1\right)}\chi\ldots\chi
x_{\left(1\cdot\left(\rho\varrho\right)\right)\left(d\right)})
\]
\[
\hskip-1pt\otimes\hskip-1pt
(x_{\left(1\cdot\left(\rho\varrho\right)\right)
\left(d+1\right)}\chi\ldots\chi
x_{\left(1\cdot\left(\rho\varrho\right)\right)
\left(d+e\right)})\hskip-2pt \otimes\ldots\otimes\hskip-2pt
(x_{\left(1\cdot\left(\rho\varrho\right)\right)
\left(p+1\right)}\chi\ldots\chi
x_{\left(1\cdot\left(\rho\varrho\right)\right) \left(n\right)})=
\]
\[
\sum_{d+e+\cdots+h=n}\sum_{\varsigma\in M\left(\chi;n;d,e,\ldots,
h \right)} \chi(\varsigma) (x_{\varsigma
\left(1\right)}\chi\ldots\chi x_{\varsigma\left(d\right)})
\]
\[
\hskip-1pt\otimes\hskip-1pt (x_{\varsigma
\left(d+1\right)}\chi\ldots\chi x_{\varsigma
\left(d+e\right)})\hskip-2pt \otimes\ldots\otimes\hskip-2pt
(x_{\varsigma \left(p+1\right)}\chi\ldots\chi
x_{\varsigma\left(n\right)})=
\]
\[
c_k(x_1\chi\ldots\chi x_n).
\]

\end{proof}

Let us denote by $m_E$ the multiplication of the algebra $[\chi](E)$:
\[
m_E\colon [\chi](E)\hskip-2pt\otimes\hskip-2pt [\chi](E)\to
[\chi](E),
\]
\[
x_1\chi\ldots\chi x_d\hskip-2pt\otimes\hskip-2pt y_1\chi\ldots\chi
y_e\mapsto x_1\chi\ldots\chi x_d\chi y_1\chi\ldots\chi y_e,
\]
and by $\varepsilon_E\colon K\to [\chi](E)$,
$\varepsilon_E(a)=a1$, the unit of the algebra $[\chi](E)$.

\begin{corollary} \label{IV.2}(i) The $K$-linear map
$c_E\colon [\chi](E)\to [\chi](E)\hskip-2pt\otimes\hskip-2pt
[\chi](E)$ defines a structure of graded coassociative
$K$-coalgebra on the graded $K$-module $[\chi](E)$, which is,
moreover, counital, with counit, the linear form $\epsilon_E$
defined by the rule
\[
\epsilon_E\colon [\chi](E)\to K,
\]
\[
\epsilon_E(z)=\left\{
\begin{array}{ll}
z & \mbox{if $z\in [\chi]^0(E)$}\\
0 & \mbox{if $z\in ([\chi](E))_+$};
 \end{array}
\right.
\]
(ii) The structure $([\chi](E),c_E,\epsilon_E)$ of graded
coassociative $K$-coalgebra with counit on the graded $K$-module
$[\chi](E)$ defines by functoriality a structure of graded
associative algebra with unit on its dual
$([\chi](E))^{*gr}=[\chi](E^*)$, and the last one coincide with
the canonical structure $([\chi](E^*),m_{E^*},\varepsilon_{E^*})$
of graded associative algebra with unit on the graded $K$-module
$[\chi](E^*)$;

\end{corollary}

\begin{proof} (i) The case $k=3$ of Proposition~\ref{IV.1} yields
coassociativity of $[\chi](E)$. We have
\[
(\epsilon_E\hskip-2pt\otimes\hskip-2pt 1)(c_E(x_1\chi\ldots\chi
x_n))=
\]
\[
(\epsilon_E\hskip-2pt\otimes\hskip-2pt
1)(\sum_{p+h=n}\sum_{\rho\in M\left(\chi;n;p,h\right)}\chi(\rho)
(x_{\rho\left(1\right)}\chi\ldots\chi
x_{\rho\left(p\right)})\hskip-2pt\otimes\hskip-2pt
(x_{\rho\left(p+1\right)}\chi\ldots\chi x_{\rho\left(n\right)}))=
\]
\[
\sum_{p+h=n}\sum_{\rho\in M\left(\chi;n;p,h\right)}\chi(\rho)
\epsilon_E((x_{\rho\left(1\right)}\chi\ldots\chi
x_{\rho\left(p\right)}))\hskip-2pt\otimes\hskip-2pt
(x_{\rho\left(p+1\right)}\chi\ldots\chi x_{\rho\left(n\right)})=
\]
\[
\sum_{\rho\in M\left(\chi;n;0,n\right)}\chi(\rho)
\epsilon_E(1)\hskip-2pt\otimes\hskip-2pt
(x_{\rho\left(1\right)}\chi\ldots\chi x_{\rho\left(n\right)})=
\]
\[
1\hskip-2pt\otimes\hskip-2pt x_1\chi\ldots\chi
x_n=x_1\chi\ldots\chi x_n.
\]
Similarly,
\[
(1\hskip-2pt\otimes\hskip-2pt\epsilon_E)(c_E(x_1\chi\ldots\chi
x_n))=
\]
\[
(1\hskip-2pt\otimes\hskip-2pt\epsilon_E)(\sum_{d+q=n}\sum_{\rho\in
M\left(\chi;n;d,q\right)}\chi(\rho)
(x_{\rho\left(1\right)}\chi\ldots\chi
x_{\rho\left(d\right)})\hskip-2pt\otimes\hskip-2pt
(x_{\rho\left(d+1\right)}\chi\ldots\chi x_{\rho\left(n\right)}))=
\]
\[
\sum_{d+q=n}\sum_{\rho\in M\left(\chi;n;d,q\right)}\chi(\rho)
(x_{\rho\left(1\right)}\chi\ldots\chi
x_{\rho\left(d\right)})\hskip-2pt\otimes\hskip-2pt
\epsilon_E((x_{\rho\left(p+1\right)}\chi\ldots\chi
x_{\rho\left(n\right)}))=
\]
\[
\sum_{\rho\in M\left(\chi;n;n,0\right)}\chi(\rho)
(x_{\rho\left(1\right)}\chi\ldots\chi
x_{\rho\left(n\right)})\hskip-2pt\otimes\hskip-2pt \epsilon_E(1)=
\]
\[
x_1\chi\ldots\chi x_n\hskip-2pt\otimes\hskip-2pt 1=
x_1\chi\ldots\chi x_n.
\]
Therefore
\[
(\epsilon_E\hskip-2pt\otimes\hskip-2pt 1)\circ
c_E=(1\hskip-2pt\otimes \hskip-2pt\epsilon_E)\circ c_E=1.
\]

(ii) Corollary~\ref{III.20} yields that the multiplication
$m_{E^*}$ in the graded algebra
$([\chi](E^*),m_{E^*},\varepsilon_{E^*})$ is the transpose of the
comultiplication $c_E$ of the graded coassociative $K$-coalgebra
with counit $([\chi](E),c_E,\epsilon_E)$. Moreover, the counit
$\epsilon_E$ is an element of $([\chi](E))^{*gr}$, such that if
$z\in [\chi](E)$, $z=z_0+z_1+z_2+\cdots$, then $\langle
z,\epsilon_E\rangle = z_0 = z_01$. The transpose of $\epsilon_E$
is the $K$-linear map $K^*\to ([\chi](E))^{*gr}$, $\ell\mapsto
\ell\circ\epsilon_E$. We compose it with the canonical isomorphism
$K\to K^*$, and, after the identification of $([\chi](E))^{*gr}$
with $[\chi](E^*)$ via the isomorphism from Theorem~\ref{III.13},
(i), we get the $K$-linear map $K\to [\chi](E^*)$, $k\mapsto k1$,
and this is the unit $1$ of the algebra $[\chi](E^*)$.

\end{proof}

\section{Inner products of a $\chi$-vector and a $\chi$-form}

\label{V}

The semi-symmetric algebra $[\chi](E)$ becomes a
$\hbox{\ccc Z}$-graded $K$-module by setting
$[\chi]^d(E)=0$ for negative integers $d$.

Let $d$ and $q\geq 0$ be integers with $d+q=n$. Let $a=a_1\chi\ldots\chi a_q$
be a fixed decomposable $q-\chi$-vector. The right multiplication by $a$ in
the algebra $[\chi](E)$,
\[
x_1\chi\ldots\chi x_d\mapsto x_1\chi\ldots\chi x_d\chi a_1\chi\ldots\chi a_q,
\]
defines an endomorphism $e'(a)$ of degree $q$ of the $\hbox{\ccc Z}$-graded $K$-module
$[\chi](E)$. The transpose of $e'(a)$ is an endomorphism $i'(a)$ of degree
$-q$ of the dual $\hbox{\ccc Z}$-graded $K$-module $[\chi](E^*)$. We define $e'(a)$
and $i'(a)$ for $a\in [\chi](E)$ by linearity.

For any $\chi$-vector $a\in [\chi](E)$ and for any $\chi$-form $a^*\in [\chi](E^*)$
denote the $\chi$-form $i'(a)(a^*)$ by $a\rfloor a^*$ and call it \emph{left inner product
of $a$ and $a^*$}. Thus,
\[
\langle x\chi a, a^*\rangle=\langle x, a\rfloor a^*\rangle
\]
for $x\in [\chi](E)$.

\begin{proposition}\label{V.1} Let $d$ and $q\geq 0$ be integers with non-negative
sum $n=d + q$. Then for any decomposable $q-\chi$-vector $a=a_1\chi\ldots\chi a_q$,
and for any decomposable $n-\chi$-form $a^*=a_1^*\chi\ldots\chi a_n^*$, the left inner
product $a\rfloor a^*$  is the $d-\chi$-linear form
\[
\sum_{\rho\in M\left(\chi;n;d,q\right)}\chi(\rho)
\langle a_1\chi\ldots\chi a_q, a_{\rho\left(d+1\right)}^*\chi\ldots\chi
a_{\rho\left(n\right)}^*\rangle a_{\rho\left(1\right)}^*\chi\ldots\chi
a_{\rho\left(d\right)}^*
\]
in case $n\geq q$, and $0$ in case $n<q$.

\end{proposition}

\begin{proof} In case $n<q$ we have $a\rfloor a^*=0$ by the definition of the endomorphism
$i'(a)$. Otherwise, $i'(a_1\chi\ldots\chi a_q)(a_1^*\chi\ldots\chi a_n^*)$ is the linear form
\[
x_1\chi\ldots\chi x_d\mapsto \langle x_1\chi\ldots\chi
x_d\chi a_1\chi\ldots\chi a_q, a_1^*\chi\ldots\chi a_n^*\rangle
\]
on $[\chi](E)$. Proposition~\ref{III.18} yields
\[
\langle x_1\chi\ldots\chi x_d\chi a_1\chi\ldots\chi a_q, a_1^*\chi\ldots\chi a_n^*\rangle=
\]
\[
\sum_{\rho\in M\left(\chi;n;d,q\right)}\chi(\rho)
\langle x_1\chi\ldots\chi x_d, a_{\rho\left(1\right)}^*\chi\ldots\chi a_{\rho\left(d\right)}^*\rangle
\langle a_1\chi\ldots\chi a_q, a_{\rho\left(d+1\right)}^*\chi\ldots\chi
a_{\rho\left(n\right)}^*\rangle=
\]
\[
\langle x_1\chi\ldots\chi x_d, \sum_{\rho\in M\left(\chi;n;d,q\right)}\chi(\rho)
\langle a_1\chi\ldots\chi a_q, a_{\rho\left(d+1\right)}^*\chi\ldots\chi
a_{\rho\left(n\right)}^*\rangle a_{\rho\left(1\right)}^*\chi\ldots\chi
a_{\rho\left(d\right)}^*\rangle.
\]
After the identification of $[\chi]^d(E))^{*gr}$ with $[\chi](E^*)$, we obtain the result.

\end{proof}

Given non-negative integers $d$, $q$ with $d+q=n$, a
integer $m\geq 1$, and $i\in J(\chi,m,d)$, $j\in J(\chi,m,q)$, $k\in J(\chi,m,n)$, one sets
\[
M_{k,.,j}(\chi;n;d,q)=\{\rho\in M(\chi;n;d,q)\mid j_1=k_{\rho\left(d+1\right)},\ldots,j_q=k_{\rho\left(n\right)}\},
\]
\[
M_{k,i,.}^\prime(\chi;n;d,q)=\{\rho\in M(\chi;n;d,q)\mid k_{\rho\left(1\right)}=i_1,\ldots,k_{\rho\left(d\right)}=i_d\}.
\]

\begin{corollary}\label{V.2}  Let $(e_\ell)_{\ell=1}^m$ be a basis for the $K$-module
$E$ and let $(e_\ell^*)_{\ell=1}^m$ be its dual basis in the dual $K$-module $E^*$.
Let $(e_j)_{j\in J\left(\chi,m\right)}$ and $(e_k^*)_{k\in J\left(\chi,m\right)}$ be the corresponding  bases
of $[\chi](E)$ and $[\chi](E^*)$, respectively. If $j\in J(\chi,m,q)$, $k\in J(\chi,m,n)$, and if $d+q=n$, then
the left inner product $e_j\rfloor e_k^*$  is the $d-\chi$-linear form
\[
\sum_{\rho\in M_{k,.,j}\left(\chi;n;d,q\right)}\chi(\rho)
e_{\rho\left(1\right)}^*\chi\ldots\chi
e_{\rho\left(d\right)}^*
\]
in case $n\geq q$, and $0$ in case $n<q$.

\end{corollary}

\begin{proof} In accord with Proposition~\ref{V.1}, in case $n<q$ we have $e_j\rfloor e_k^*=0$,
and in case $n\geq q$, we have
\[
e_j\rfloor e_k^*=
\]
\[
\sum_{\rho\in M\left(\chi;n;d,q\right)}\chi(\rho)
\langle e_{j_1}\chi\ldots\chi e_{j_q}, e_{k_{\rho\left(d+1\right)}}^*\chi\ldots\chi
e_{k_{\rho\left(n\right)}}^*\rangle e_{k_{\rho\left(1\right)}}^*\chi\ldots\chi
e_{k_{\rho\left(d\right)}}^*=
\]
\[
\sum_{\rho\in M_{k,.,j}\left(\chi;n;d,q\right)}\chi(\rho)
e_{k_{\rho\left(1\right)}}^*\chi\ldots\chi
e_{k_{\rho\left(d\right)}}^*.
\]

\end{proof}

\begin{proposition}\label{V.4} The addition and the external composition law
$(a, a^*)\mapsto a\rfloor a^*$ on $[\chi](E^*)$ define on this set a structure of
left unital $[\chi](E)$-module.

\end{proposition}

\begin{proof} The external composition law is  bilinear and the
associativity of the the graded algebra $[\chi](E)$ is equivalent to  the
equality $e'(a\chi b)=e'(b)\circ e'(a)$ for $a, b\in [\chi](E)$.  Then  $i'(a\chi b)=
i'(a)\circ i'(b)$, and hence $(a\chi b)\rfloor a^*=a\rfloor (b\rfloor a^*)$. Moreover,
$1\rfloor a^*=a^*$.

\end{proof}

Let $p\geq 0$ and $h$ be integers with $p+h=n$. Let $a^*=a_1^*\chi\ldots\chi a_p^*$
be a fixed decomposable $p-\chi$-form. The left multiplication by $a^*$ in
the algebra $[\chi](E^*)$,
\[
x_1^*\chi\ldots\chi x_h^*\mapsto a_1^*\chi\ldots\chi a_p^*\chi x_1^*\chi\ldots\chi x_h^*,
\]
defines an endomorphism $e(a^*)$ of degree $p$ of the $\hbox{\ccc Z}$-graded $K$-module
$[\chi](E^*)$. The transpose of $e(a^*)$ is an endomorphism $i(a^*)$ of degree
$-p$ of the $\hbox{\ccc Z}$-graded $K$-module $[\chi](E)$. We define $e(a^*)$
and $i(a^*)$ for $a^*\in [\chi](E)$ by linearity.

For any $\chi$-form $a^*\in [\chi](E^*)$, and for any $\chi$-vector $a\in [\chi](E)$
denote the $\chi$-vector $i(a^*)(a)$ by $a\lfloor a^*$ and call it \emph{right inner product
of $a$ and $a^*$}. Thus,
\[
\langle a\lfloor a^*, x^*\rangle=\langle a, a^* \chi x^* \rangle
\]
for $x^*\in [\chi](E^*)$.

\begin{proposition}\label{V.5} Let $h$ and $p\geq 0$ be integers with non-negative
sum $n=p+h$. Then for any decomposable $n-\chi$-vector $a=a_1\chi\ldots\chi a_n$,
and for any decomposable $p-\chi$-form $a^*=a_1^*\chi\ldots\chi a_p^*$, the right inner
product $a\lfloor a^*$  is the $h-\chi$-vector
\[
\sum_{\rho\in M\left(\chi;n;p,h\right)}\chi(\rho)
\langle a_{\rho\left(1\right)}\chi\ldots\chi a_{\rho\left(p\right)},
a_1^*\chi\ldots\chi a_p^*\rangle a_{\rho\left(p+1\right)}\chi\ldots\chi a_{\rho\left(n\right)}
\]
in case $n\geq p$, and $0$ in case $n<p$.

\end{proposition}

\begin{proof} In case $n<p$ we have $a\lfloor a^*=0$ by the definition of the endomorphism
$i(a^*)$. Otherwise, according to Proposition~\ref{III.18} we have
\[
\langle a\lfloor a^*, x_1^*\chi\ldots\chi x_h^*\rangle=
\langle a_1\chi\ldots\chi a_n, a_1^*\chi\ldots\chi a_p^*\chi x_1^*\chi\ldots\chi x_h^*\rangle=
\]
\[
\sum_{\rho\in M\left(\chi;n;p,h\right)}\chi(\rho)
\langle a_{\rho\left(1\right)}\chi\ldots\chi a_{\rho\left(p\right)}, a_1^*\chi\ldots\chi a_p^*\rangle
\langle a_{\rho\left(p+1\right)}\chi\ldots\chi a_{\rho\left(n\right)}, x_1^*\chi\ldots\chi x_h^*\rangle=
\]
\[
\langle \sum_{\rho\in M\left(\chi;n;p,h\right)}\chi(\rho)
\langle a_{\rho\left(1\right)}\chi\ldots\chi a_{\rho\left(p\right)},
a_1^*\chi\ldots\chi a_p^*\rangle a_{\rho\left(p+1\right)}\chi\ldots\chi
a_{\rho\left(n\right)}, x_1^*\chi\ldots\chi x_h^*\rangle,
\]
and we get the result.

\end{proof}

\begin{corollary}\label{V.6}  Let $(e_\ell)_{\ell=1}^m$ be a basis for the $K$-module
$E$ and let $(e_\ell^*)_{\ell=1}^m$ be its dual basis in the dual $K$-module $E^*$.
Let $(e_j)_{j\in J\left(\chi,m\right)}$ and $(e_k^*)_{k\in J\left(\chi,m\right)}$ be the corresponding  bases
of $[\chi](E)$ and $[\chi](E^*)$, respectively. If $j\in J(\chi,m,p)$, $k\in J(\chi,m,n)$, and if $p+h=n$, then
the right inner product $e_k\lfloor e_j^*$  is the $h-\chi$-vector
\[
\sum_{\rho\in M_{k,j,.}\left(\chi;n;p,h\right)}\chi(\rho)
e_{k_{\rho\left(p+1\right)}}\chi\ldots\chi e_{k_{\rho\left(n\right)}}
\]
in case $n\geq p$, and $0$ in case $n<p$.

\end{corollary}

\begin{proof} In accord with Proposition~\ref{V.1}, in case $n<p$ we have $e_k\lfloor e_j^*=0$,
and in case $n\geq p$, we have
\[
e_k\rfloor e_j^*=
\]
\[
\sum_{\rho\in M\left(\chi;n;p,h\right)}\chi(\rho)
\langle e_{k_{\rho\left(1\right)}}\chi\ldots\chi
e_{k_{\rho\left(p\right)}}, e_{j_1}^*\chi\ldots\chi e_{j_p}^*\rangle e_{k_{\rho\left(p+1\right)}}\chi\ldots\chi
e_{k_{\rho\left(n\right)}}=
\]
\[
\sum_{\rho\in M_{k,j,.}\left(\chi;n;p,h\right)}\chi(\rho)
e_{k_{\rho\left(p+1\right)}}\chi\ldots\chi e_{k_{\rho\left(n\right)}}.
\]

\end{proof}

\begin{proposition}\label{V.7} The addition and the external composition law
$(a, a^*)\mapsto a\lfloor a^*$ on $[\chi](E)$ define on this set a structure of
right unital $[\chi](E^*)$-module.

\end{proposition}

\begin{proof} The external composition law is  bilinear and the
associativity of the the graded algebra $[\chi](E^*)$ is equivalent to  the
equality $e(a^*\chi b^*)=e(a^*)\circ e(b^*)$ for $a^*, b^*\in [\chi](E^*)$.  Then  $i(a^*\chi b^*)=
i(b^*)\circ i(a^*)$, and hence $a\lfloor (a^*\chi b^*)=(a\lfloor a^*)\lfloor b^*$. Moreover,
$a\lfloor 1=a$.

\end{proof}

 \section{Acknowledgements}
 This work was supported in part by Grant MM-1503/2005 of the 
Bulgarian Foundation of Scientific Research.

\appendix

\section{Appendix}

\label{A}

\begin{notation}\label{A.2} Let $d$, $e$,$\ldots$, $h$ be $k$ in number nonnegative
integers with $d+e+\cdots +h=n$. We assume $k\leq n$. Let
$\alpha\colon [1,d]\to [1,n]$, $\beta\colon [1,e]\to
[1,n]$,$\ldots$, $\gamma\colon [1,h]\to [1,n]$, be strictly
increasing maps with disjoint images. Let $\theta_\alpha\in S_n$
be a permutation with $\theta_\alpha(1)=\alpha(1)$,$\ldots$,
$\theta_\alpha(d)=\alpha(d)$, let $\theta_\beta\in S_n$ be a
permutation with $\theta_\beta(1)=\beta(1)$,$\ldots$,
$\theta_\beta(e)=\beta(e)$,$\ldots$, let $\theta_\gamma\in S_n$ be
a permutation with $\theta_\gamma(1)=\gamma(1)$,$\ldots$,
$\theta_\gamma(h)=\gamma(h)$. For any permutation $\theta\in S_n$
we denote by $c_\theta\colon S_n\to S_n$ the conjugation
$c_\theta(\zeta)=\theta\zeta\theta^{-1}$. We have
\[
c_{\theta_\alpha}(S_d)=S_{I\hskip-1pt m\alpha},\mbox{\ }
c_{\theta_\beta}(S_e)=S_{I\hskip-1pt m\beta},\mbox{\ }\ldots,
c_{\theta_\gamma}(S_h)=S_{I\hskip-1pt m\gamma}.
\]
Let $K$ be a commutative ring with unit $1$. Let $U\leq S_d$,
$V\leq S_e$,$\ldots$ ,$W\leq S_h$ be permutation groups, and let
$\varepsilon\colon U\to U(K)$, $\delta\colon V\to U(K)$,$\ldots$,
$\varpi\colon W\to U(K)$, be linear $K$-valued characters. We
embed the Cartesian product $U\times V\times\cdots\times W$ in
$S_n$ as $X=c_{\theta_\alpha}(U)c_{\theta_\beta}(V)\ldots
c_{\theta_\gamma}(W)$ and for any $\zeta\in X$,
$\zeta=c_{\theta_\alpha}(\sigma)c_{\theta_\beta}(\tau)\ldots 
c_{\theta_\gamma}(\eta)$, $\sigma\in U$, $\tau\in V$,$\ldots$,
$\eta\in W$, we set
\[
\chi(\zeta)=\varepsilon(\sigma)\delta(\tau)\ldots\varpi(\eta).
\]
The map $\chi\colon X\to U(K)$ is a $K$-linear character of the
group $X$. Let $E$ be a $K$-module and let $(x_1,\ldots,x_d)\in
E^d$, $(y_1,\ldots,y_e)\in E^e$,$\ldots$, $(z_1,\ldots,z_h)\in
E^h$ be generic elements. We set
\[
\xi_i= \left\{
\begin{array}{llll}
x_{\alpha^{-1}\left(i\right)} & \mbox{if $i\in I\hskip-1pt m\alpha$}\\
y_{\beta^{-1}\left(i\right)} & \mbox{if $i\in I\hskip-1pt m\beta$}\\
\vdots & \vdots\\
z_{\gamma^{-1}\left(i\right)} & \mbox{if $i\in I\hskip-1pt m\gamma$}\\
 \end{array}
\right.
\]

Let $Y\leq S_n$ be a permutation group with $X\leq
Y$, and let $M_{U,V,\ldots,W}^{\alpha,\beta,\ldots,\gamma}
\hskip-2pt\left(Y\right)$ be the set of all
lexicographically minimal representatives of the left classes of
$Y$ modulo $X$. For any $\zeta'\in Y$, $\zeta\in
M_{U,V,\ldots,W}^{\alpha,\beta,\ldots,\gamma}
\hskip-2pt\left(Y\right)$, we denote by $\zeta'\cdot\zeta$
the lexicographically minimal representative of $\zeta'\zeta$
modulo $X$, and set $\zeta'\cdot\zeta=\zeta'\zeta\upsilon_{\zeta'\zeta}$,
where $\upsilon_{\zeta'\zeta}\in X$,
$\upsilon_{\zeta'\zeta}=c_{\theta_\alpha}(\sigma)c_{\theta_\beta}(\tau)\ldots
c_{\theta_\gamma}(\eta)$, with $\sigma\in U$, $\tau\in
V$,$\ldots$, $\eta\in W$.

In case an $\omega$-invariant sequence of characters
$\chi=(\chi_d)_{d\geq 1}$ is given, if the opposite is not stated,
we specialize the maps $\alpha,\beta,\ldots,\gamma$, the groups
$U,V,\ldots, W$, and the characters $\varepsilon, \delta,\ldots,
\varpi$, on them, as follows: $\alpha(1)=1$,$\ldots$,
$\alpha(d)=d$, $\beta(1)=d+1$,$\ldots$, $\beta(e)=d+e$,$\ldots$,
$\gamma(1)=d+e+\cdots+1$,$\ldots$, $\gamma(h)=d+e+\cdots+h$,
$U=W_d$, $V=W_e$, $\ldots$, $W=W_h$, $Y=W_n$,
$\varepsilon=\chi_d$, $\delta=\chi_e$,$\ldots$, $\varpi=\chi_h$.
Then
\[
c_{\theta_\alpha}(U)=W_d,\mbox{\ }
c_{\theta_\beta}(V)=\omega^d(W_e),\mbox{\ }\ldots,
c_{\theta_\gamma}(W)=\omega^{d+e+\cdots}(W_h),
\]
and, using notation from Remark~\ref{III.16},
\[
M_{U,V,\ldots,W}^{\alpha,\beta,\ldots,\gamma}
\hskip-2pt\left(Y\right)=M(\chi;n;d,e,\ldots,h).
\]

\end{notation}

\begin{lemma} \label{A.3} The rule $(\zeta',\zeta)\mapsto \zeta'\cdot\zeta$
defines a left action of the group $Y$ on the set
$M\hskip-2pt\left(Y;\alpha,\beta,\ldots,\gamma\right)$.

\end{lemma}

\begin{proof} Let $\zeta''\in Y$. The three elements
$(\zeta''\zeta')\cdot\zeta$, $\zeta''(\zeta'\cdot\zeta)$, and
$\zeta''\cdot(\zeta'\cdot\zeta)$ are in the class
$\zeta''\zeta'\zeta X$, so we get $(\zeta''\zeta')\cdot\zeta=
\zeta''\cdot(\zeta'\cdot\zeta)$. Finally,
$1_{Y}\cdot\zeta=\zeta$.

\end{proof}

\begin{lemma} \label{A.4} Let $\pi$ be a linear $K$-valued character of $Y$, and
$\pi_{\mid X}=\chi$. Let $\varepsilon^2=1_U$,
$\delta^2=1_V$,$\ldots$, $\varpi^2=1_W$, and
$\pi^2=1_{Y}$. The formula
\[
[\pi]^n(E)\to\coprod_{d+e+\cdots+h=n}[\varepsilon]^d(E)\hskip-2pt\otimes\hskip-2pt
[\delta]^e(E)\ldots\otimes\hskip-2pt [\varpi]^h(E),
\]
\[
\xi_1\pi\ldots\pi\xi_n\mapsto
\sum_{d+e+\cdots+h=n}\sum_{\zeta\in
M_{U,V,\ldots,W}^{\alpha,\beta,\ldots,\gamma}
\hskip-2pt\left(Y\right)}\pi(\zeta)
\]
\[ (\xi_{\zeta\left(\alpha\left(1\right)\right)}\varepsilon\ldots\varepsilon
\xi_{\zeta\left(\alpha\left(d\right)\right)})\hskip-2pt\otimes\hskip-2pt
(\xi_{\zeta\left(\beta\left(1\right)\right)}\delta\ldots\delta
\xi_{\zeta\left(\beta\left(e\right)\right)})\hskip-2pt\otimes\ldots\otimes\hskip-2pt
(\xi_{\zeta\left(\gamma\left(1\right)\right)}\varpi\ldots\varpi
\xi_{\zeta\left(\gamma\left(h\right)\right)}),
\]
defines a $K$-linear map.

\end{lemma}

\begin{proof} The map
\[
f\colon
E^n\to\coprod_{d+e+\cdots+h=n}[\varepsilon]^d(E)\hskip-2pt\otimes\hskip-2pt
[\delta]^e(E)\ldots\otimes\hskip-2pt [\varpi]^h(E),
\]
\[
f(\xi_1,\ldots,\xi_n)=\sum_{d+e+\cdots+h=n}\sum_{\zeta\in
M_{U,V,\ldots,W}^{\alpha,\beta,\ldots,\gamma}
\hskip-2pt\left(Y\right)}\pi(\zeta)
\]
\[ (\xi_{\zeta\left(\alpha\left(1\right)\right)}\varepsilon\ldots\varepsilon
\xi_{\zeta\left(\alpha\left(d\right)\right)})\hskip-2pt\otimes\hskip-2pt
(\xi_{\zeta\left(\beta\left(1\right)\right)}\delta\ldots\delta
\xi_{\zeta\left(\beta\left(e\right)\right)})\hskip-2pt\otimes\ldots\otimes\hskip-2pt
(\xi_{\zeta\left(\gamma\left(1\right)\right)}\varpi\ldots\varpi
\xi_{\zeta\left(\gamma\left(h\right)\right)}),
\]
is multilinear and semi-symmetric of weight $\pi$. Indeed,
let $\zeta'\in Y$. We have
\[
f(\xi_{\zeta'\left(1\right)},\ldots,\xi_{\zeta'\left(n\right)})=
\sum_{d+e+\cdots+h=n}\sum_{\zeta\in
M_{U,V,\ldots,W}^{\alpha,\beta,\ldots,\gamma}
\hskip-2pt\left(Y\right)}\pi(\zeta)
\]
\[ (\xi_{\zeta'\left(\zeta\left(\alpha\left(1\right)\right)\right)}\varepsilon\ldots\varepsilon
\xi_{\zeta'\left(\zeta\left(\alpha\left(d\right)\right)\right)})\hskip-2pt\otimes\hskip-2pt
(\xi_{\zeta'\left(\zeta\left(\beta\left(1\right)\right)\right)}\delta\ldots\delta
\xi_{\zeta'\left(\zeta\left(\beta\left(e\right)\right)\right)})\hskip-2pt\otimes\ldots
\]
\[
\otimes\hskip-2pt
(\xi_{\zeta'\left(\zeta\left(\gamma\left(1\right)\right)\right)}
\varpi\ldots\varpi\xi_{\zeta'\left(\zeta\left(\gamma\left(h\right)\right)\right)})=
\]
\[
\pi(\zeta')\sum_{d+e+\cdots+h=n}\sum_{\zeta\in
M_{U,V,\ldots,W}^{\alpha,\beta,\ldots,\gamma}
\hskip-2pt\left(Y\right)}\pi(\zeta'\zeta)
\]
\[ (\xi_{\zeta'\left(\zeta\left(\alpha\left(1\right)\right)\right)}\varepsilon\ldots\varepsilon
\xi_{\zeta'\left(\zeta\left(\alpha\left(d\right)\right)\right)})\hskip-2pt\otimes\hskip-2pt
(\xi_{\zeta'\left(\zeta\left(\beta\left(1\right)\right)\right)}\delta\ldots\delta
\xi_{\zeta'\left(\zeta\left(\beta\left(e\right)\right)\right)})\hskip-2pt\otimes\ldots
\]
\[
\otimes\hskip-2pt
(\xi_{\zeta'\left(\zeta\left(\gamma\left(1\right)\right)\right)}
\varpi\ldots\varpi\xi_{\zeta'\left(\zeta\left(\gamma\left(h\right)\right)\right)}).
\]
Since
\[
\pi(\zeta'\cdot\zeta)=\pi(\zeta'\zeta\upsilon_{\zeta'\zeta})=
\pi(\zeta'\zeta)\chi(\upsilon_{\zeta'\zeta})=
\pi(\zeta'\zeta)\varepsilon(\sigma)\delta(\tau)\ldots
\varpi(\eta),
\]
using Lemma~\ref{A.3}, we have
\[
f(\xi_{\zeta'\left(1\right)},\ldots,\xi_{\zeta'\left(n\right)})=
\pi(\zeta')\sum_{d+e+\cdots+h=n}\sum_{\zeta\in
M_{U,V,\ldots,W}^{\alpha,\beta,\ldots,\gamma}
\hskip-2pt\left(Y\right)}\pi(\zeta'\zeta)
\]
\[ (\xi_{\zeta'\left(\zeta\left(\alpha\left(1\right)\right)\right)}\varepsilon\ldots\varepsilon
\xi_{\zeta'\left(\zeta\left(\alpha\left(d\right)\right)\right)})\hskip-2pt\otimes\hskip-2pt
(\xi_{\zeta'\left(\zeta\left(\beta\left(1\right)\right)\right)}\delta\ldots\delta
\xi_{\zeta'\left(\zeta\left(\beta\left(e\right)\right)\right)})\hskip-2pt\otimes\ldots
\]
\[
\otimes\hskip-2pt
(\xi_{\zeta'\left(\zeta\left(\gamma\left(1\right)\right)\right)}
\varpi\ldots\varpi\xi_{\zeta'\left(\zeta\left(\gamma\left(h\right)\right)\right)})=
\]
\[
\pi(\zeta')\sum_{d+e+\cdots+h=n}\sum_{\zeta\in
M_{U,V,\ldots,W}^{\alpha,\beta,\ldots,\gamma}
\hskip-2pt\left(Y\right)}\pi(\zeta'\cdot\zeta)
\]
\[
\varepsilon(\sigma)
(\xi_{\zeta'\left(\zeta\left(\alpha\left(1\right)\right)\right)}\varepsilon\ldots\varepsilon
\xi_{\zeta'\left(\zeta\left(\alpha\left(d\right)\right)\right)})\hskip-2pt\otimes\hskip-2pt
\delta(\tau)(\xi_{\zeta'\left(\zeta\left(\beta\left(1\right)\right)\right)}\delta\ldots\delta
\xi_{\zeta'\left(\zeta\left(\beta\left(e\right)\right)\right)})\hskip-2pt\otimes\ldots
\]
\[
\otimes\hskip-2pt
\varpi(\eta)(\xi_{\zeta'\left(\zeta\left(\gamma\left(1\right)\right)\right)}
\varpi\ldots\varpi\xi_{\zeta'\left(\zeta\left(\gamma\left(h\right)\right)\right)})=
\]
\[
\pi(\zeta')\sum_{d+e+\cdots+h=n}\sum_{\zeta\in
M_{U,V,\ldots,W}^{\alpha,\beta,\ldots,\gamma}
\hskip-2pt\left(Y\right)}\pi(\zeta'\cdot\zeta)
\]
\[
(\xi_{\zeta'\left(\zeta\left(\alpha\left(\sigma\left(1\right)\right)\right)\right)}\varepsilon\ldots\varepsilon
\xi_{\zeta'\left(\zeta\left(\alpha\left(\sigma\left(d\right)\right)\right)\right)})\hskip-2pt\otimes\hskip-2pt
(\xi_{\zeta'\left(\zeta\left(\beta\left(\tau\left(1\right)\right)\right)\right)}\delta\ldots\delta
\xi_{\zeta'\left(\zeta\left(\beta\left(\tau\left(e\right)\right)\right)\right)})\hskip-2pt\otimes\ldots
\]
\[
\otimes\hskip-2pt
\varpi(\eta)(\xi_{\zeta'\left(\zeta\left(\gamma\left(\eta\left(1\right)\right)\right)\right)}
\varpi\ldots\varpi\xi_{\zeta'\left(\zeta\left(\gamma\left(\eta\left(h\right)\right)\right)\right)})=
\]
\[
\pi(\zeta')\sum_{d+e+\cdots+h=n}\sum_{\zeta\in
M_{U,V,\ldots,W}^{\alpha,\beta,\ldots,\gamma}
\hskip-2pt\left(Y\right)}\pi(\zeta'\cdot\zeta)
\]
\[
(\xi_{\zeta'\left(\zeta\left(\upsilon_{\zeta'\zeta}\left(\alpha\left(1\right)\right)\right)\right)}
\varepsilon\ldots\varepsilon
\xi_{\zeta'\left(\zeta\left(\upsilon_{\zeta'\zeta}\left(\alpha\left(d\right)\right)\right)\right)})
\hskip-2pt\otimes\hskip-2pt
(\xi_{\zeta'\left(\zeta\left(\upsilon_{\zeta'\zeta}\left(\beta\left(1\right)\right)\right)\right)}
\delta\ldots\delta
\xi_{\zeta'\left(\zeta\left(\upsilon_{\zeta'\zeta}\left(\beta\left(e\right)\right)\right)\right)})
\hskip-2pt\otimes\ldots
\]
\[
\otimes\hskip-2pt
(\xi_{\zeta'\left(\zeta\left(\upsilon_{\zeta'\zeta}\left(\gamma\left(1\right)\right)\right)\right)}
\varpi\ldots\varpi\xi_{\zeta'\left(\zeta\left(\upsilon_{\zeta'\zeta}\left(\gamma\left(h\right)\right)
\right)\right)})=
\]
\[
\pi(\zeta')\sum_{d+e+\cdots+h=n}\sum_{\zeta\in
M_{U,V,\ldots,W}^{\alpha,\beta,\ldots,\gamma}
\hskip-2pt\left(Y\right)}\pi(\zeta'\cdot\zeta)
\]
\[
(\xi_{\left(\zeta'\cdot\zeta\right)\left(\alpha\left(1\right)\right)}\varepsilon\ldots\varepsilon
\xi_{\left(\zeta'\cdot\zeta\right)\left(\alpha\left(d\right)\right)})\hskip-2pt\otimes\hskip-2pt
(\xi_{\left(\zeta'\cdot\zeta\right)\left(\beta\left(1\right)\right)}
\delta\ldots\delta
\xi_{\left(\zeta'\cdot\zeta\right)\left(\beta\left(e\right)\right)})\hskip-2pt\otimes\ldots
\]
\[
\otimes\hskip-2pt
(\xi_{\left(\zeta'\cdot\zeta\right)\left(\gamma\left(1\right)\right)}
\varpi\ldots\varpi\xi_{\left(\zeta'\cdot\zeta\right)\left(\gamma\left(h\right)\right)})=
\]
\[
\pi(\zeta')f(\xi_1,\ldots,\xi_n).
\]
Therefore, according to \cite[(1.1.1)]{[3]}, $f$ gives rise to the desired
$K$-linear map.

\end{proof}

Let $\chi=(\chi_d)_{d\geq 1}$ be an $\omega$-invariant sequence of
characters. Using Notation~\ref{A.2}, we have

\begin{lemma} \label{A.5} The maps
\[
M\left(\chi;n;p,h\right)\times M\left(\chi;p;d, e,\ldots\right)\to
M\left(\chi;n;d, e,\ldots, h\right),
\]
\[
M\left(\chi;n;d,q\right)\times \omega^d M\left(\chi;q;e,\ldots,
h\right)\to M\left(\chi;n;d, e,\ldots, h\right),
\]
\[
(\rho,\varrho)\mapsto 1\cdot(\rho\varrho),
\]
are bijections.

\end{lemma}

\begin{proof} If $W_n/W_p\times\omega^p(W_h)$ is a set of
representatives of the left classes of $W_n$ modulo
$W_p\times\omega^p(W_h)$, if
$W_p\times\omega^p(W_h)/W_d\times\omega^d(W_e)\times\cdots\times\omega^p(W_h)$
is a set of representatives of the left classes of
$W_p\times\omega^p(W_h)$ modulo
$W_d\times\omega^d(W_e)\times\cdots\times\omega^p(W_h)$, then the
family
\[
\{\rho\varrho\mid (\rho,\varrho)\in
(W_n/W_p\times\omega^p(W_h))\times
(W_p\times\omega^p(W_h)/W_d\times\omega^d(W_e)\times\cdots\times\omega^p(W_h))\}
\]
of elements of $W_n$ is a set of representatives of the left
classes of $W_n$ modulo
$W_d\times\omega^d(W_e)\times\cdots\times\omega^p(W_h)$. Thus, the
first map is a bijection because $M\left(\chi;p;d,
e,\ldots\right)$ is a set of representatives of the left classes
of $W_p\times\omega^p(W_h)$ modulo
$W_d\times\omega^d(W_e)\times\cdots\times\omega^p(W_h)$.
Similarly, if $W_n/W_d\times\omega^d(W_q)$ is a set of
representatives of the left classes of $W_n$ modulo
$W_d\times\omega^d(W_q)$, if
$W_d\times\omega^d(W_q)/W_d\times\omega^d(W_e)\times\cdots\times\omega^p(W_h)$
is a set of representatives of the left classes of
$W_d\times\omega^d(W_q)$ modulo
$W_d\times\omega^d(W_e)\times\cdots\times\omega^p(W_h)$, then the
family
\[
\{\rho\varrho\mid (\rho,\varrho)\in
(W_n/W_d\times\omega^d(W_q))\times
(W_d\times\omega^d(W_q)/W_d\times\omega^d(W_e)\times\cdots\times\omega^p(W_h))\}
\]
of elements of $W_n$ is a set of representatives of the left
classes of $W_n$ modulo
$W_d\times\omega^d(W_e)\times\cdots\times\omega^p(W_h)$. The
second map is a bijection, too, because $\omega^d
M\left(\chi;q;e,\ldots, h\right)$ is a set of representatives of
the left classes of $W_d\times\omega^d(W_q)$ modulo
$W_d\times\omega^d(W_e)\times\cdots\times\omega^p(W_h)$.

\end{proof}

\end{document}